\RequirePackage{ifpdf}
\ifpdf 
\documentclass[pdftex]{sigma}
\else
\documentclass{sigma}
\fi

\renewcommand\t{\tilde}
\newcommand\cn{\mbox{cn}}
\newcommand\dn{\mbox{dn}}
\newcommand\sn{\mbox{sn}}

\numberwithin{equation}{section}

\begin{document}

\allowdisplaybreaks

\renewcommand{\thefootnote}{$\star$}

\renewcommand{\PaperNumber}{033}

\FirstPageHeading

\ShortArticleName{Frobenius Determinant and Biorthogonal Polynomials}

\ArticleName{Elliptic Hypergeometric Laurent Biorthogonal\\
Polynomials with a Dense Point Spectrum\\ on the Unit Circle\footnote{This paper is a contribution to the Proceedings of the Workshop ``Elliptic Integrable Systems, Isomonodromy Problems, and Hypergeometric Functions'' (July 21--25, 2008, MPIM, Bonn, Germany). The full collection
is available at
\href{http://www.emis.de/journals/SIGMA/Elliptic-Integrable-Systems.html}{http://www.emis.de/journals/SIGMA/Elliptic-Integrable-Systems.html}}}

\Author{Satoshi TSUJIMOTO~$^\dag$ and Alexei ZHEDANOV~$^\ddag$}

\AuthorNameForHeading{S.~Tsujimoto and A.~Zhedanov}

\Address{$^\dag$~Department of Applied Mathematics and Physics,
Graduate School of Informatics,\\
\hphantom{$^\dag$}{}~Kyoto University, Kyoto 606-8501, Japan}
\EmailD{\href{mailto:tsujimoto@amp.i.kyoto-u.ac.jp}{tsujimoto@amp.i.kyoto-u.ac.jp}}

\Address{$^\ddag$~Donetsk Institute for Physics and Technology,
Donetsk 83114, Ukraine}
\EmailD{\href{mailto:zhedanov@kinetic.ac.donetsk.ua}{zhedanov@kinetic.ac.donetsk.ua}}

\ArticleDates{Received November 30, 2008, in f\/inal form March 15,
2009; Published online March 19, 2009}

\Abstract{Using the technique of the elliptic Frobenius
determinant, we construct new elliptic solutions of the
$QD$-algorithm. These solutions can be interpreted as elliptic
solutions of the discrete-time Toda chain as well. As a
by-product, we obtain new explicit orthogonal and biorthogonal
polynomials in terms of the elliptic hypergeometric function
${_3}E_2(z)$. Their recurrence coef\/f\/icients are expressed in terms
of the elliptic functions. In the degenerate case we obtain the
Krall--Jacobi polynomials and their biorthogonal analogs.}


\Classification{33E05; 33E30; 33C47}

\renewcommand{\thefootnote}{\arabic{footnote}}
\setcounter{footnote}{0}

\section{Introduction}

 In this paper we present new explicit
solutions for the two-point $QD$-algorithm \cite{Cabe} (which is
equivalent to the discrete-time relativistic Toda chain \cite{Ru,Suris,KMZ}). These solutions can be naturally
constructed starting from the famous Frobenius elliptic
determinant (see, e.g.,~\cite{FrSt,Amde}). This approach
allows one to f\/ind an explicit expression for corresponding
Laurent biorthogonal polynomials in terms of the elliptic
hypergeometric function ${_3}E_2(z)$. These polynomials contain
several free parameters and appear to be biorthogonal on the unit
circle with respect to a~dense point spectrum. In two special
cases we already obtained explicit examples of $\cn$- and $\dn$-elliptic polynomials which are orthogonal on the unit circle with
respect to a positive dense point measure \cite{Zhe_cndn}. These
polynomials provide f\/irst known explicit (i.e.\ expressed in terms
of the elliptic hypergeometric function) examples of such measures
(see also~\cite{Simon} for general properties of polynomials
orthogonal with respect to measures of singular type and~\cite{Mag} for an example of such polynomials). The obtained
polynomials $P_n(z)$ possess a remarkable ``classical'' property.
This means that ${\cal D} P_n(z) = \mu_n \t P_{n-1}(z),$ where
$\cal D$ is a generalized derivative operator: ${\cal D} z^n=
\mu_n z^{n-1}$ (with some coef\/f\/icients $\mu_n$) and $\t P_n(z)$
are polynomials of the same type but with shifted parameters. In
our case the operator $\cal D$ is an elliptic generalization of
the ordinary derivative operator ${\cal D} = \partial_z$, $\mu_n
=n$ and $q$-derivative operator with $\mu_n = (q^n-1)/(q-1)$.

In the degenerated case (when both periods of elliptic functions
become inf\/inity) we obtain biorthogonal analogs of the
Krall--Jacobi orthogonal polynomials. We show that these
biorthogo\-nal polynomials satisfy a 4th order dif\/ferential
equation which can be presented in the form of quadratic operator
pencil.

\section{Laurent biorthogonal polynomials and their basic
properties}

 The Laurent biorthogonal
polynomials LBP $P_n(z)$ appeared in problems connected with the
two-points Pad\'e approximations (see, e.g.,~\cite{JT}).

We shall recall their def\/inition and general properties (see, e.g.,
\cite{JT,HR,IsMas}, where equivalent Laurent {\it
orthogonal functions} are considered).

Let $\cal L$ be some linear functional def\/ined on all possible
monomials $z^n$ by the moments
\begin{gather*} c_n = {\cal L} \{z^n\}, \qquad n
=0 , \pm 1, \pm 2 \dots .
\end{gather*}
In general the moments
$c_n$ are arbitrary complex numbers. The functional $\cal L$ is
def\/ined on the space of Laurent polynomials ${\cal P}(z)=
\sum\limits_{n=-N_1}^{N_2}a_n z^n$ where $a_n$ are arbitrary complex
numbers and~$N_{1,2}$ arbitrary integers:
\[
{\cal L}\{{\cal P}(z)\}= \sum_{n=-N_1}^{N_2}a_nc_n .
\]

The monic LBP $P_n(z)$ are def\/ined by the determinant \cite{HR}
\begin{gather}
P_n(z)=(\Delta_n)^{-1} \left | \begin{array}{cccc}
c_0 & c_1 & \dots & c_n \\ c_{-1} & c_0 & \dots & c_{n-1} \\ \dots
& \dots & \dots & \dots\\ c_{1-n}& c_{2-n}& \dots & c_1\\ 1& z &
\dots & z^n  \end{array} \right |, \label{deterP}
\end{gather} where
$\Delta_n$ is def\/ined as the Toeplitz determinant
\[
 \Delta_n = \left |
\begin{array}{cccc} c_0 & c_{1} & \dots &
c_{n-1}\\ c_{-1}& c_0 & \dots & c_{n-2}\\ \dots & \dots & \dots & \dots\\
c_{1-n} & c_{2-n} & \dots & c_0 \end{array} \right | .
\] It is obvious from def\/inition~\eqref{deterP} that the
polynomials $P_n(z)$ satisfy the orthogonality property
\begin{gather*} {\cal
L}\{P_n(z) z^{-k}\} =h_n \delta_{kn}, \qquad 0 \le k \le n,
\end{gather*}
where the normalization constants $h_n$ are \begin{gather*}
h_0=c_0, \qquad h_n = \Delta_{n+1}/\Delta_{n}. 
 \end{gather*} This
orthogonality property can be rewritten as the biorthogonal
relation \cite{Pas, HR},
\begin{gather*}
{\cal L}\{P_n(z) Q_m(1/z)\} =
h_n \delta_{nm}, 
\end{gather*} where the polynomials $Q_n(z)$ are
def\/ined by the formula
\begin{gather}  Q_n(z)=(\Delta_n)^{-1} \left |
\begin{array}{cccc} c_0 & c_{-1} & \dots & c_{-n} \\ c_{1} & c_0 &
\dots & c_{1-n} \\ \dots & \dots & \dots & \dots\\ c_{n-1}&
c_{n-2}& \dots & c_{-1}\\ 1& z & \dots & z^n  \end{array} \right
|. \label{deterQ}
\end{gather} We note that the polynomials $Q_n(z)$ are
again LBP with moments $c^{\{Q\}}_n = c_{-n}$.

In what follows we will assume that \begin{gather} \Delta_n \ne 0, \qquad
n=1,2,\dots \label{inDe} \end{gather} and that \begin{gather} \Delta_n^{(1)} \ne 0,
\qquad n=1,2,\dots , \label{inT} \end{gather} where by $\Delta_n^{(j)}$ we
denote the determinants
\begin{gather}
 \Delta_0^{(j)}=1, \qquad
 \Delta_n^{(j)} = \left | \begin{array}{cccc} c_j & c_{j+1} &
\dots & c_{n+j-1}\\ c_{j-1}& c_{j} & \dots & c_{n+j-2}\\ \dots &
\dots & \dots & \dots\\ c_{1+j-n} & c_{2+j-n} & \dots & c_j
\end{array} \right | .\label{Dj}
\end{gather}

If the conditions \eqref{inDe} and \eqref{inT} are fulf\/illed, the
polynomials $P_n(z)$ satisfy the recurrence relation (see, e.g.,~\cite{HR})
\begin{gather} P_{n+1}(z) +(d_n-z)P_n(z) =zb_nP_{n-1}(z), \qquad n
\ge 1, \label{recP} \end{gather} where the recurrence coef\/f\/icients are
\begin{gather}
 d_n = -\frac{P_{n+1}(0)}{P_n(0)} =h_n^{-1} \frac{T_{n+1}}{T_n}
=\frac{T_{n+1} \Delta_n}{T_n \Delta_{n+1}} \ne 0, \qquad n=0,1,\dots, \label{d_n} \\
 b_n= d_n \frac{h_n}{h_{n-1}} = \frac{T_{n+1} \Delta_{n-1}}{T_n
\Delta_{n}} \ne 0, \qquad n=1,2,\dots \label{b_n}
\end{gather} with $T_n=
\Delta_n^{(1)}$. Note the important relation \begin{gather*} \frac{b_n}{d_n} =
\frac{h_n}{h_{n-1}} = \frac{\Delta_{n-1}\Delta_{n+1}}{\Delta_n^2},
\qquad n=1,2,\dots 
\end{gather*} from which one can obtain
expression for the normalization constant $h_n$ in terms of the
recurrence parameters: \begin{gather} h_n = \prod_{i=1}^n \frac{b_i}{d_i}.
\label{hbd} \end{gather}

There is a one-to-one correspondence between the moments $c_n$ and
the recurrence coef\/f\/i\-cients~$b_n$,~$d_n$ (provided restrictions~$b_n d_n \ne 0$ are fulf\/illed).

We say that the LBP are regular if $b_n d_n \ne 0$. This condition
is equivalent to the condition \begin{gather*} \Delta_n \Delta^{(1)}_n \ne 0,
\qquad n=0,1,\dots . \label{reg_Delta}
\end{gather*}

In the regular case there is a simple formula relating the
biorthogonal partners $Q_n(z)$ with polynomials $P_n(z)$
\cite{HR}: \begin{gather} Q_n(z) = \frac{zP_{n+1}(1/z) -
z^{n-1}P_n(1/z)}{P_n(0)}. \label{Q_from_P} \end{gather}

In what follows we will use so-called rescaled LBP \begin{gather*} \t P_n(z) =
q^n P_n(z/q), \qquad n=0,1,\dots 
\end{gather*} with some non-zero
parameter $q$. It is easily verif\/ied that the rescaled polynomials
$\t P_n(z)$ are monic LBP satisfying the recurrence relation
\[
\t P_{n+1}(z) +(\t d_n-z) \t P_n(z) =z \t b_n \t P_{n-1}(z)
\]
with \begin{gather*} \t d_n = q b_n, \qquad \t b_n = q b_n. 
\end{gather*} The
rescaled LBP $\t P_n(z)$ dif\/fer from initial LBP $P_n(z)$ only by
a trivial rescaling of recurrence parameters. The moments $\t c_n$
of the rescaled LBP are connected with initial moments $c_n$ by
the relation $\t c_n = q^n c_n$. Note that the rescaled
biorthogonal partners $Q_n(z)$ are transformed as \begin{gather} \t Q_n(z) =
q^{-n} Q_n(zq). \label{res_Q} \end{gather}

There is a connection between the LBP and the restricted
relativistic Toda chain~\cite{KMZ}. Assume that LBP $P_n(z;t)$
depend on an additional (so-called ``time'') parameter $t$. This
mean that the recurrence coef\/f\/icients $b_n(t)$, $d_n(t)$ become
functions of the parameter $t$. We assume that the relation \begin{gather*}
\dot P_n(z) = -\frac{b_n}{d_n} P_{n-1}(z) 
\end{gather*} holds for
all $n=0,1,\dots$. This ansatz leads to the following equations
for the recurrence coef\/f\/icients~\cite{KMZ}
\begin{gather}
 \dot d_n =  \frac{b_{n+1}}{d_{n+1}} - \frac{b_n}{d_{n-1}}, \qquad
 \dot b_n = b_n\left( \frac{1}{d_{n}}  -
\frac{1}{d_{n-1}}\right). \label{1rte}
\end{gather} For the corresponding
moments $c_n(t)$ we have the relation \begin{gather*} \dot c_n = c_{n-1},
\qquad n=0, \pm 1, \pm2, \dots . 
\end{gather*}

Another possible ansatz \cite{KMZ} \begin{gather*} \dot P_n(z) =
-b_n(P_n(z)-zP_{n-1}(z))
\end{gather*} leads to the equations
\begin{gather}
 \dot d_n= -d_n (b_{n+1} - b_n), \qquad
 \dot b_n = -b_n(b_{n+1}-b_{n-1}+d_{n-1}-d_n). \label{2rte}
 \end{gather} In
this case we have for the moments the relation \begin{gather*} \dot c_n =
c_{n+1}, \qquad n=0, \pm 1, \pm 2, \dots . 
\end{gather*}

In spite of the apparent dif\/ference between equations \eqref{1rte}
and \eqref{2rte}, it can be shown (see, e.g.,~\cite{KMZ}) that these
two systems are both equivalent to the restricted relativistic
Toda chain equations. The term ``restricted'' in this context means
that it is assumed an additional condition \begin{gather*} b_0 =0. 
\end{gather*} This means that in formulas \eqref{1rte} or \eqref{2rte} we should
assume $n=0,1,2,\dots$. For the nonrestricted relativistic Toda
chain equations \eqref{1rte} or \eqref{2rte} are valid for all integer
values of $n=0, \pm 1, \pm 2, \dots$.

\section[Laurent biorthogonal polynomials and $QD$-algorithm]{Laurent biorthogonal polynomials and $\boldsymbol{QD}$-algorithm}

 The (restricted) ``discrete-time''
relativistic Toda chain corresponds to the following ansatz for
the moments \begin{gather*} c_n(t+h) = c_{n+1}(t), \qquad n=0, \pm 1, \pm 2,
\dots, 
\end{gather*} where $h$ is an arbitrary parameter. We
have the transformation formula for the corresponding Laurent
biorthogonal polynomials \begin{gather} P_n(z;t+h) = P_n(z;t) + b_n(t)
P_{n-1}(z;t) \label{dt_P} \end{gather} and \begin{gather} (d_n - b_n)   P_n(z;t-h) =
zP_n(z;t) - P_{n+1}(z;t). \label{dt_P2} \end{gather} Formulas \eqref{dt_P} and
\eqref{dt_P2} can be interpreted as Christof\/fel and Geronimus
transformations for LBP \cite{Zh_LBP}.

The corresponding recurrence coef\/f\/icients are transformed as
\cite{Zh_LBP} \begin{gather} d_n(t+h) = d_{n-1}
\frac{b_{n+1}-d_n}{b_n-d_{n-1}}, \qquad b_n(t+h) = b_{n}
\frac{b_{n+1}-d_n}{b_n-d_{n-1}} \label{dt_bd} \end{gather} (in r.h.s.\ of
\eqref{dt_bd} it is assumed the argument $t$ for the coef\/f\/icients
$b_n$, $d_n$). These relations can be presented in a slightly
dif\/ferent equivalent form as \begin{gather} b_n \tilde d_n = d_{n-1} \t b_n,
\qquad \t b_n - \t d_n = b_{n+1} - d_n,  \label{rel_QD} \end{gather} where we
have denoted $\t b_n = b_n(t+h)$ etc for brevity. Relations
\eqref{rel_QD} describe so-called $QD$-algorithm for the two-point
Pad\'e approximation (see, e.g.,~\cite{Cabe} for details). In other
words, the (restricted) discrete-time relativistic Toda chain is
equivalent to the $QD$-algorithm for the two-point Pad\'e
approximation.

Usually, this algorithm works as follows. We start from the given
moments $c_n(t)$, $n=0, \pm 1, \pm 2\, \dots$ where the
dependence on ``time'' is trivial: \begin{gather*} c_n(t+h) = c_{n+1}(t)
\label{c:linear}
\end{gather*} and def\/ine the coef\/f\/icient $d_0(t)$ for all $t
=t_0 + jh$, $j =0 ,\pm 1, \pm 2$ as
\[ d_0(t)
=\frac{c_0(t+h)}{c_0(t)}.
\]
The initial value $t_0$ is not essential, usually it is assumed
that $t_0=0$, in this case we can write
\[
d_0(t+jh) \equiv d_0^{(j)} = \frac{c_{j+1}}{c_j}.
\]
Assume that $b_0(t)=0$ for all $t$. Then at the f\/irst step we f\/ind
$b_1(t)=b_1^{(j)}$ for all $t=jh$ from the second relation
\eqref{rel_QD}:
\[
b_1^{(j)} = d_0^{(j)} - d_0^{(j+1)}.
\]
Then we f\/ind $d_1^{(j)}$ from the f\/irst relation \eqref{rel_QD}
\[
d_1^{(j+1)} = \frac{b_1^{(j+1)} d_0^{(j)}}{b_1^{(j)}}.
\]
This process can be continued to f\/ind $b_2^{(j)}$, $d_2^{(j)}$,
$\dots$. The process is non-degenerate if $b_n^{(j)}d_n^{(j)} \ne
0$ for all $n$ and $j$. Then we obtain all sequences $d_n^{(j)}$,
$b_n^{(j)}$, $n=0,1,2,\dots$ for $j =0, \pm 1, \pm 2, \dots$.

There is a remarkable connection with the $QD$-algorithm for the
ordinary orthogonal polynomials \cite{Cabe}. Indeed, let us
introduce the monic polynomials \begin{gather} W_n^{(j)}(z) \equiv
P_n^{(j+n)}(z), \label{W_P_cor} \end{gather} where the polynomials
$P_n^{(j)}(z)$ are def\/ined as $P_n^{(j)}(z)=P_n(z;hj)$.

Then relations \eqref{dt_P} and \eqref{dt_P2} become \begin{gather} W_n^{(j-1)}(z)
= W_n^{(j)}(z) - f_n^{(j)}W_{n-1}^{(j)}(z) \label{Q_qd1} \end{gather} and \begin{gather}
zW_n^{(j+1)}(z) = W_{n+1}^{(j)}(z) - e_n^{(j+1)} W_n^{(j)}(z),
\label{Q_qd2} \end{gather} where
\[
f_n^{(j)}= b_n^{(j+n-1)}, \qquad e_n^{(j)} = b_n^{(j+n)} -
d_n^{(j+n)}.
\]
These relations can be interpreted as Geronimus and Christof\/fel
transforms for the orthogonal polynomials $W_n^{(j)}(z)$. The
compatibility condition between~\eqref{Q_qd1} and~\eqref{Q_qd2} leads to
the recurrence relation \begin{gather*} W_{n+1}^{(j)}(z) + g_n^{(j)}
W_{n}^{(j)}(z) + u_n^{(j)} W_{n-1}^{(j)}(z) = zW_{n}^{(j)}(z),
\end{gather*} which describes the three-term recurrence
relation for the ordinary orthogonal polyno\-mials $W_{n}^{(j)}(z)$
where the recurrence coef\/f\/icients are~\cite{Cabe} \begin{gather*}
g_n^{(j)}=-e_n^{(j)} -f_{n+1}^{(j)}, \qquad u_n^{(j)} =
e_n^{(j)}f_n^{(j)}. 
\end{gather*} Moreover we have compatibility
conditions for the coef\/f\/icients $e_n^{(j)}$, $f_n^{(j)}$ \begin{gather}
e_{n-1}^{(j+1)}f_n^{(j+1)} = e_n^{(j)}f_n^{(j)}, \qquad e_n^{(j+1)}
+ f_n^{(j+1)} = e_n^{(j)} + f_{n+1}^{(j)}. \label{QDR} \end{gather} Relations
\eqref{QDR} coincide with those introduced by Rutishauser and
describing the ordinary $QD$-algorithm~\cite{FF}. It is easy to
verify that relations~\eqref{QDR} are equivalent to relations~\eqref{rel_QD} for the two-point $QD$-algorithm.

Thus starting from known solution $P_n^{(j)}(z)$, $b_n^{(j)}$,
$d_n^{(j)}$ of the discrete-time relativistic Toda chain (or,
equivalently, two-point $QD$-algorithm) we can obtain a set of the
ordinary orthogonal polynomials $W_n^{(j)}(z)$ depending on
additional ``time'' parameter~$j$. Note that sometimes the
introduced orthogonal polynomials~$W_n^{(j)}(z)$ depending on an
additional discrete parameter $j$ are called the Hadamard
polynomials \cite{AC,Henrici}\footnote{The authors are indebted to
A.~Magnus for drawing their attention to these references.}.

From the def\/inition \eqref{W_P_cor} it follows that the orthogonal
polynomials $W_n^{(j)}(z)$ can be presented in determinantal form
as \begin{gather*}
 W_n^{(j)}(z)= \frac{1}{H_n^{(j)}} \left |
\begin{array}{cccc} c_{j+1} & c_{j+2} & \dots & c_{n+j+1} \\ c_{j+2} & c_{j+3} &
\dots & c_{j+n+2} \\ \dots & \dots & \dots & \dots\\ c_{n+j}&
c_{n+j+1}& \dots & c_{2n+j}\\ 1& z & \dots & z^n  \end{array}
\right |, 
\end{gather*}
where $H_n^{(j)}$ stands for the Hankel
determinant
\begin{gather*} H_n^{(j)} = \left |
\begin{array}{cccc} c_{j+1} & c_{j+2} & \dots &
c_{j+n}\\ c_{j+2}& c_{j+3} & \dots & c_{j+n+1}\\ \dots & \dots & \dots & \dots\\
c_{j+n} & c_{j+n+1} & \dots & c_{j+2n-1} \end{array} \right |
.
\end{gather*}
 Clearly we have the relation \begin{gather*} H_n^{(j)} =
(-1)^{n(n-1)/2} \Delta_n^{(n+j)}. 
\end{gather*} Thus the
orthogonal polynomials $W_n^{(j)}(z)$ are orthogonal \begin{gather*} \langle
\tau^{(j)}, W_n^{(j)}(z) W_m^{(j)}(z) \rangle = q_n^{(j)}
\delta_{nm}, 
\end{gather*} where the linear functional
$\tau^{(j)}$ is def\/ined by the moments \begin{gather*} \tau_n^{(j)} \equiv
\langle \tau^{(j)}, z^n \rangle = c_{n+j+1}, \qquad n=0,1,2,\dots, \qquad
j=0, \pm 1, \pm 2, \dots . 
\end{gather*} The normalization
constant $q_n^{(j)}$ has the expression
\[
q_n^{(j)} = \frac{H_{n+1}^{(j)}}{H_n^{(j)}} = (-1)^n
\frac{\Delta_{n+1}^{(j+n+1)}}{\Delta_n^{(j+n)}}.
\]

It would be instructive to interpret (\ref{dt_P}) and
(\ref{dt_P2}) in terms of so-called bilinear technique by using
the determinantal identities. This technique is standard in the
theory of integrable systems.

As a f\/irst step, we give a compressed expression to $d_n - b_n$ as
\begin{gather*}
 d_n -b_n = \frac{\Delta_{n+1}^{(1)}\Delta_{n}^{(-1)}}
                  {\Delta_{n+1}^{(0)}\Delta_{n}^{(0)}},
\end{gather*} which can be derived from the determinantal identity, or
Jacobi identity, for the Toeplitz determinant:
\begin{gather*}
 \Delta_{n+1}^{(j)} \Delta_{n-1}^{(j)}
 = (\Delta_n^{(j)})^2 -  \Delta_{n}^{(j+1)}\Delta_{n}^{(j-1)}.
\end{gather*}
Then the relations (\ref{dt_P})  and (\ref{dt_P2}) can be
transformed to the following bilinear equations,
\begin{gather*}
 \Delta_n^{(j)} \sigma_n^{(j+1)}
 =
 \Delta_n^{(j+1)} \sigma_n^{(j)}
 -  \Delta_{n+1}^{(j+1)}\sigma_{n-1}^{(j+1)}  ,\\
   \sigma_{n+1}^{(j)} \Delta_n^{(j)}
= z \sigma_{n}^{(j)} \Delta_{n+1}^{(j)}
 - \Delta_{n+1}^{(j+1)} \sigma_{n}^{(j-1)},
\end{gather*}
respectively, where the functions $\sigma_n^{(j)}$ are def\/ined by
\begin{gather*}
 \sigma_n^{(j)} =\left |
\begin{array}{cccc} c_0^{(j)} & c_{1}^{(j)} & \dots & c_{n}^{(j)} \\ c_{-1}^{(j)} & c_0^{(j)} &
\dots & c_{n-1}^{(j)} \\ \dots & \dots & \dots & \dots\\
c_{1-n}^{(j)}& c_{2-n}^{(j)}& \dots & c_{1}^{(j)}\\ 1& z & \dots &
z^n  \end{array} \right |. 
\end{gather*} (Note that
$\sigma_n^{(j)}$ is proportional to the Laurent biorthogonal
polynomial $P_n^{(j)}(z)$.)

\section{Laurent and Baxter biorthogonal polynomials}

There is an alternative (but essentially
equivalent) approach to biorthogonal polynomials proposed by
G.~Baxter~\cite{Baxter}. The pair $P_n(z)$, $Q_n(z)$ of the
biorthogonal polynomials is def\/ined in this approach by means of
initial conditions $P_0=Q_0=1$ and the following recurrence system
\begin{gather} P_{n+1}(z) = zP_n(z) - e_n^{(1)} Q^*_n(z), \qquad  Q_{n+1}(z) =
zQ_n(z) - e_n^{(2)} P^*_n(z),  \label{Bax_rec} \end{gather} where
$e_n^{(1,2)}$ are some complex coef\/f\/icients. It is clear that
$e_n^{(1)}=-P_{n+1}(0)$, $e_n^{(2)}=-Q_{n+1}(0)$. Notation
$P_n^*(z)$ is standard for so-called reciprocal polynomials, i.e.\
$P_n^*(z)=z^nP_n(1/z)$, $Q_n^*(z)=z^nQ_n(1/z)$. Assume that
$e_n^{(1)}e_n^{(2)}(1-e_n^{(1)}e_n^{(2)}) \ne 0$ (this is the
nondegenerate case). Then, excluding $Q_n^*(z)$ from the system
\eqref{Bax_rec} we arrive at a 3-term recurrence relation for the
polyno\-mials~$P_n(z)$:
\[
P_{n+1}(z) + d_n P_n(z) = z(P_n(z) + b_n P_{n-1}(z))
\]
coinciding with \eqref{recP}, where
\[
d_n = - \frac{e_n^{(1)}}{e_{n-1}^{(1)}}, \qquad b_n =
-\frac{e_n^{(1)}}{e_{n-1}^{(1)}}(1-{e_{n-1}^{(1)}}{e_{n-1}^{(2)}}).
\]
Clearly, polynomials $Q_n(z)$ satisfy similar relations with
interchanging superscripts~$1$,~$2$.

Conversely, assume that we have the nondegenerate Laurent
biorthogonal polynomials $P_n(z)$ satisfying \eqref{recP}. We can
construct their biorthogonal partners $Q_n(z)$ by \eqref{Q_from_P}.
Then it is elementary to verify that polynomials $P_n(z)$, $Q_n(z)$
satisfy system~\eqref{Bax_rec} with $e_n^{(1)}=-P_{n+1}(0)$, $e_n^{(2)}=-Q_{n+1}(0)$. Sometimes system \eqref{Bax_rec} is more
convenient for analysis due to apparent symmetry between
polynomials $P_n(z)$, $Q_n(z)$ and corresponding coef\/f\/icients
$e_n^{(1)}$, $e_n^{(2)}$. Note also that the Laurent and Baxter
biorthogonal polynomials in turn are equivalent to the so-called
Laurent {\it orthogonal} polynomials proposed by Jones and Thron~\cite{JT}. The Jones and Thron polynomials contains terms~$z^k$
with both positive and negative degree $k$. For details of this
equivalence see, e.g.,~\cite{HR} and~\cite{Pas}.

There is an important special case when all the Toeplitz
determinants are positive $\Delta_n >0$ and moreover the moments
satisfy the condition \begin{gather*} \bar c_n = c_{-n} 
\end{gather*}
(as usual, $\bar c_n$ means complex conjugation of $c_n$). In this
case the biorthogonal partners $Q_n(z)$ coincide with complex
conjugated polynomials $Q_n(z) = \bar P_n(z)$ and there exists
nondecreasing function $\sigma(\theta)$ of bounded variation on
the unit circle such that the orthogonality relation \begin{gather}
\int_{0}^{2 \pi} P_n\big(e^{i\theta}\big) \bar P_m(e^{-i\theta}) d
\sigma(\theta) = h_n   \delta_{nm} \label{ort_UC} \end{gather} holds. I.e.\
in this case we have polynomials $P_n(z)$ which are {\it
orthogonal} on the unit circle (abbreviated as OPUC~\cite{Simon}).
Historically, these polynomials were introduced f\/irst by Szeg\H{o}
\cite{Sz} and are called the Szeg\H{o} polynomials orthogonal on
the unit circle. They satisfy the recurrence relation \begin{gather}
P_{n+1}(z) = zP_n(z) - a_n z^n \bar P_n(1/z), \label{Sz_rec} \end{gather}
where the coef\/f\/icients $a_n = -P_{n+1}(0)$ are called the
ref\/lection (or Schur, or Verblunsky, \dots) parameters. The relation
\eqref{Sz_rec} was f\/irst derived by Szeg\H{o} himself~\cite{Sz}. The
ref\/lection parameters are complex numbers satisfying the important
inequality \begin{gather} |a_n|<1, \qquad n=0,1,2,\dots . \label{a_n_res} \end{gather} In
fact, condition \eqref{a_n_res} is equivalent to the condition of
positive def\/inite Toeplitz forms $\Delta_n >0$ or to existence of
a positive measure on the unit circle providing orthogonality
property~\eqref{ort_UC}.

If, additionally, all the moments are real, then they satisfy
condition $c_{-n}=c_n$. In this case the ref\/lection parameters are
real parameters satisfying the restriction $-1 < a_n <1$, $
n=0,1,2,\dots$. The biorthogonal partners then coincide with
initial polynomials $Q_n(z) = P_n(z)$. It is easy to show that the
measure $d\sigma$ is symmetric with respect to real axis in this
case, namely the function~$\sigma(\theta)$ satisf\/ies the condition
$\sigma(2\pi -\theta) + \sigma(\theta) = {\rm const}$.

For further details concerning theory of OPUC see, e.g.,~\cite{Ger,Simon}.

\section{Frobenius elliptic determinant formula\\ and biorthogonal functions}

 Assume that $v_i$, $u_i$, $i=0,1,\dots$ are
two arbitrary sequences of complex numbers. Let \begin{gather*}
H_n=\det||g_{ij} ||_{i,j=0,\dots,n-1}, 
\end{gather*} where \begin{gather*} g_{ij}=
\frac{\sigma(u_i+v_j+\beta)}{\sigma(u_i+v_j)\sigma(\beta)}
\exp(\gamma_1 u_i + \gamma_2 v_j), 
\end{gather*} where
$\sigma(z)$ is the standard Weierstrass sigma function (see, e.g.,~\cite{Akhiezer,WW}) and $\beta$, $\gamma_1$, $\gamma_2$ are
arbitrary.

Recall that the Weierstrass sigma function is def\/ined by the
inf\/inite product \cite{Akhiezer} \begin{gather*} \sigma(u)= \Pi'
\left(1-\frac{u}{s}\right)  \exp\left(\frac{u}{s} + \frac{u^2}{2
s^2}\right), 
\end{gather*} where the product is taken over
all points of the lattice $s= 2m \omega_1 + 2 m' \omega_3$, $
m,m'=0, \pm 1, \pm 2, \dots$ excluding the point with $m=m'=0$. $2
\omega_1$ and $2 \omega_3$ are the so-called primitive elliptic
periods. It is convenient to introduce the third period $2
\omega_2 = -2 \omega_1 - 2 \omega_3$~\cite{Akhiezer}. The
Weierstrass sigma function possess quasi-periodic properties
\cite{Akhiezer} \begin{gather*} \sigma(u+ 2 \omega_{\alpha}) = - \exp(2
\eta_{\alpha}(u+ \omega_{\alpha}) ) \sigma(u), \qquad
\alpha=1,2,3,
\end{gather*} where the constants
$\eta_{\alpha}$ are def\/ined as
\[
\eta_{\alpha} = \zeta(\omega_{\alpha}), \qquad \alpha=1,2,3
\]
and $\zeta(u)= \sigma'(u)/\sigma(u)$ is the Weierstrass zeta
function \cite{Akhiezer}.

We have \begin{gather} H_n = \frac{\sigma(U+V+\beta)
\prod\limits_{i>j}{\sigma(u_i-u_j)\sigma(v_i-v_j)}}{\sigma(\beta)
\prod\limits_{i,j}\sigma(u_i+v_j)}   \exp(\gamma_1 U + \gamma_2 V)
\label{det_H} \end{gather} (we denote $U=\sum\limits_{i=0}^{n-1}u_i$,
$V=\sum\limits_{i=0}^{n-1}v_i$ for simplicity, moreover it is assumed that
the upper limit for $i,j$ in the products is $n-1$).

Formula \eqref{det_H} was obtained by Frobenius in \cite{Fr}. A
simple elementary method to derive formu\-la~\eqref{det_H} can be found
in \cite{Amde}. Frobenius and Stickelberger derived also in
\cite{FrSt} several other explicit formulas for ``elliptic
determinants'' in connection with the theory of rational
interpolation.

Let $\phi_k(x)$, $\psi_k(x)$, $k=0,1,\dots$ (we assume that
$\phi_0= \psi_0=1$) be two sets of functions in some argument $x$.
Assume that there exists a linear functional $\cal L$ such that
\begin{gather*} \langle{\cal L}, \phi_j(x) \psi_i(x)\rangle = g_{ij}.
 \end{gather*}

The linear functional $\cal L$ is def\/ined on the
space of functions constructed from bilinear combinations of the
type
\[
f(x) = \sum_{i,k=0} c_{ik} \phi_i(x) \psi_k(x)
\]
with arbitrary coef\/f\/icients $c_{ik}$.

Introduce the following functions \begin{gather}
P_n(x)=\frac{1}{\Delta_n}\left |\begin{array}{cccc} g_{00} &
g_{01} & \dots &
g_{0n}\\ g_{10}& g_{11} & \dots & g_{1n}\\ \dots & \dots & \dots & \dots\\
g_{n-1,0} & g_{n-1,1} & \dots & g_{n-1,n} \\ \phi_0(x) & \phi_1(x)
& \dots & \phi_n(x)
\end{array} \right | \label{P_det}
\end{gather} and \begin{gather*}
 Q_n(x)=\frac{1}{\Delta_n}\left |\begin{array}{cccc}
g_{00} & g_{10} & \dots &
g_{n0}\\ g_{01}& g_{11} & \dots & g_{n1}\\ \dots & \dots & \dots & \dots\\
g_{0,n-1} & g_{1,n-1} & \dots & g_{n,n-1} \\ \psi_0(x) & \psi_1(x)
& \dots & \psi_n(x)
\end{array} \right |, 
\end{gather*} where \begin{gather} \Delta_n=H_n=\det||g_{ij} ||_{i,j=0,\dots,n-1}. \label{DH}
\end{gather} By construction, these functions are biorthogonal \begin{gather*} \langle
{\cal L}, P_n(x) Q_m(x) \rangle = h_n \delta_{nm} 
\end{gather*} with respect to the functional $\cal L$, where the
normalization coef\/f\/icients $h_n$ are
\[
h_n = \frac{\Delta_{n+1}}{\Delta_{n}}.
\]
Expanding the determinant in \eqref{P_det} over the last row we have
explicit expression for the polynomial $P_n(x)$: \begin{gather*} P_n(x) =
\sum_{k=0}^n (-1)^{n-k} p_{nk} \phi_k(x), 
\end{gather*} where
\begin{gather*} p_{nk}  = \frac{H_n(k)}{\Delta_n}. 
\end{gather*} The auxiliary
determinants $H_n(k)$ are def\/ined by canceling the $k$th column,
i.e. \begin{gather*} H_n(k) = \det||g_{ij}(k)||_{i,j=0,\dots,n-1}, \label{Hnk}
\end{gather*}
where \begin{gather*} g_{ij}(k)
=\frac{\sigma(u_i+v_j(k)+\beta)}{\sigma(u_i+v_j(k))\sigma(\beta)}
\exp(\gamma_1 u_i + \gamma_2 v_j(k)). 
\end{gather*} Here the
sequence $v_i(k)$ is def\/ined as \begin{gather*} v_i(k)=\left\{\begin{array}{ll}  v_{i} \  &
\mbox{if} \ \ i<k,   \\ v_{i+1} \  &  \mbox{if} \ \  i \ge k.
\end{array} \right . 
\end{gather*}
Thus the determinant $H_n(k)$ is obtained
from the determinant $H_n$ by replacing sequence $v_i$ with the
sequence $v_i(k)$. (By def\/inition $H_n(n) = H_n$ and $v_i(n) =
v_i$.) But formula \eqref{det_H} is valid for {\it any} sequences
$u_i$, $v_i$. Hence we can calculate all the determinant $H_n(k)$
explicitly. Omitting obvious calculations we present the result
\begin{gather*} p_{nk} = e^{\gamma_2(v_n-v_k)}
\frac{\sigma(U+V+v_n-v_k+\beta)}{\sigma(U+V+\beta)}
\left[{n\atop k }\right]
\prod\limits_{i=0}^{n-1}\frac{\sigma(u_i+v_k)}{\sigma(u_i+v_n)}   ,
\end{gather*}
where \begin{gather*} \left[{n\atop k }\right] =
\frac{\prod\limits_{i=0}^{n-1}\sigma(v_n-v_i)}{\prod_{i=0}^{k-1}\sigma(v_k-v_i)
  \prod\limits_{i=k+1}^{n}\sigma(v_i-v_k) } 
  \end{gather*} are
``generalized binomial coef\/f\/icients''. Similar expression can be
obtained for the biorthogonal partners $Q_n(x)$ if one replaces
the parameters $v_i$ with $u_i$.

In case when the sequence
$v_j$ is {\it linear} with respect to $j$: $v_j=w j + \xi$ we
obtain the conventional ``elliptic binomial coef\/f\/icients''
\cite{GR}: \begin{gather*} \left[{n\atop k }\right] = \frac{[n]!}{[k]![n-k]!}
= (-1)^k \frac{[-n]_k}{[1]_k}, 
\end{gather*}
 where $[x] = \sigma(w x)/\sigma(w)$ is so-called ``elliptic number'' and
$[x]_k =[x] [x+1] \cdots [x+k-1]$ is elliptic Pochhammer symbol.
Note that usually the elliptic number is def\/ined in terms of the
theta function $[x] = \theta_1(w x)/\theta_1(w)$ \cite{GR}, but
for our purposes these def\/initions are in fact equivalent.

We thus constructed an explicit system of biorthogonal functions
$P_n(x)$, $Q_n(x)$ starting from the elliptic Frobenius determinant.
This system can be further specif\/ied by a concrete choice of the
basic functions $\phi_n(x)$, $\psi_n(x)$ and the linear functional
$\sigma$. Note that the idea to construct explicit families of
biorthogonal functions directly from corresponding Gram
determinants is due to Wilson~\cite{Wilson}. For general
biorthogonal rational functions the determinant representation can
be found e.g.\ in \cite{SZ} and \cite{DZ}.

\section{Laurent biorthogonal polynomials\\ from the Frobenius determinant}

In what follows we will assume that the
period $2 \omega_1$ is a real while the period $2 \omega_3$ is
purely imagi\-nary. This means that the fundamental parallelogram is
a rectangle. Such choice is standard for many practical purposes
because in this case the function $\sigma(x)$ takes real values on
the real axis~$x$~\cite{Akhiezer}. This is important for existence
of a positive orthogonality measure on the unit circle.

Put
\[
\gamma_1=\gamma_2=\gamma, \qquad u_i = -iw + \alpha, \qquad v_j = jw,
\]
where $w$ is an arbitrary real parameter which is incommensurable
with the real period $2 \omega_1$ over the integers, i.e.\ we will
assume that \begin{gather} w N_1 \ne \omega_1 N_2 \label{incomp_w} \end{gather} for any
integers $N_1$, $N_2$. Then for the entries of the Frobenius matrix
we have \begin{gather*} g_{ij} = e^{\gamma w(j-i) + \gamma \alpha}
\frac{\sigma(w(j-i) + \beta + \alpha) }{\sigma(w(j-i) + \alpha)
\sigma(\beta)}. 
\end{gather*} This matrix has the Toeplitz form.
We can therefore def\/ine corresponding monic Laurent biorthogonal
polynomials by the formula \begin{gather} P_n(z)=\frac{1}{\Delta_n}\left |
\begin{array}{cccc} c_0 & c_1 & \dots &
c_n\\ c_{-1}& c_0 & \dots & c_{n-1}\\ \dots & \dots & \dots & \dots\\
c_{-n+1} & c_{-n+2} & \dots & c_{1} \\ 1 & z & \dots & z^n
\end{array} \right |, \label{PL_det}
\end{gather} where the moments are def\/ined as \begin{gather} c_n = g_{0,n} = e^{\gamma
wn + \gamma \alpha}  \frac{\sigma(wn + \beta + \alpha)
}{\sigma(wn + \alpha) \sigma(\beta)} \label{mom_T} \end{gather}
 and the Toeplitz determinant $\Delta_n$ is def\/ined by~\eqref{DH}.

As in the previous section, def\/ine the elliptic numbers $[x]$ as
\[
[x] = \sigma(wx)/\sigma(w),
\]
and the elliptic Pochhammer symbol
\[
[x]_n =[x][x+1] \cdots[x+n-1].
\]
The elliptic hypergeometric function is def\/ined by the formula \begin{gather}
{_{r+1}}E_r\left ({{\vec a} \atop {\vec b}}; z \right) =
\sum_{s=0}^{\infty} \frac{[a_1]_s [a_2]_s \cdots [a_{r+1}]_s}{[1]_s
[b_1]_s [b_2]_s \cdots [b_r]_s}   e^{Ms(s-1)}  z^s, \label{Ge} \end{gather}
where
\[
M=\frac{\eta_1}{2 \omega_1} w^2 \left(1 + \sum_{i=1}^r b_i - \sum_{i=1}^{r+1} a_i\right).
\]
We have
\begin{proposition}
The Laurent biorthogonal polynomials defined by formulas
\eqref{PL_det} and \eqref{mom_T} are expressed in terms of the elliptic
hypergeometric function: \begin{gather} P_n(z) = B_n \, {_{3}}E_2 \left( {-n
, \hat \alpha+1, -(\hat \alpha+1) n - \hat \beta +1 \atop \hat
\alpha+1 - n, -(\hat \alpha+1) n - \hat \beta}; ze^{-\gamma   w}
\right), \label{hyp_LP} \end{gather} where $\hat \alpha=\alpha w^{-1}$, $\hat \beta=\beta w^{-1}$ and
\begin{gather} B_n = e^{\gamma wn}
\frac{[-\hat \alpha]_n}{[\hat \alpha +1]_n}   \frac{[\hat \alpha
n + \hat \beta + n]}{[\hat \alpha n + \hat \beta]} \label{B-cf} \end{gather}
is the coefficient to provide monicity $P_n(z) = z^n + O(z^{n-1})$
of the polynomials $P_n(z)$.
\end{proposition}

\noindent
{\bf Remark.} The parameters of the elliptic hypergeometric
function in our case satisfy condition
\[
1 + b_1 + b_2 = a_1 + a_2 + a_3
\]
and hence $M=0$ in the def\/inition of the hypergeometric function~\eqref{Ge}. Our def\/inition of the elliptic hypergeometric function is
in accordance with the conventional one~\cite{GR,Spi_Es}. The main dif\/ference is replacing the theta
functions with the Weierstrass sigma functions. This replacement
leads to appearance of the additional factor $e^{Ms(s-1)}$.
Indeed, there is relation between these functions \cite{Akhiezer}
\[
\sigma(z) = {\rm const} \cdot \exp\left(\frac{\eta_1 z^2}{2\omega_1}\right)
\theta_1(z/(2\omega_1))
\]
(the constant factor is not essential because it is canceled in
all expressions for elliptic hypergeometric series). Using this
relation we can replace all sigma functions with the theta
functions~$\theta_1(z)$ which leads to formula~\eqref{Ge}.

\medskip

Now we calculate the normalization coef\/f\/icients $h_n$ directly
from Frobenius formula~\eqref{det_H}: \begin{gather} h_n =
\frac{\Delta_{n+1}}{\Delta_n} = \frac{e^{\gamma \alpha}
}{\sigma(\alpha)}   \frac{\sigma(\alpha(n+1) +
\beta)}{\sigma(\alpha   n + \beta)}   \frac{[n]!^2}{[-\hat
\alpha+1]_n [\hat \alpha+1]_n}. \label{h_expl} \end{gather} In what follows
we will assume the following restriction $\alpha \ne wm$ for any
integers $m$. Indeed, otherwise the normalization coef\/f\/icient
$h_n$ becomes singular and we have a degeneration.

We observe also that the determinant $\Delta_n^{(1)}$ def\/ined by
\eqref{Dj} with $j=1$ is obtained from $\Delta_n^{(1)}$ by the shift
of the parameter $\alpha \to \alpha+w$ because $c_{n+1}(\alpha) =
c_n(\alpha+w)$. Thus in general we have the important formula \begin{gather*}
\Delta_n^{(j)}(\alpha) = \Delta_n(\alpha+jw). 
 \end{gather*}
In particular, we have \begin{gather} h_n^{(1)} =
\frac{\Delta_{n+1}^{(1)}}{\Delta_{n}^{(1)}} = \frac{T_{n+1}}{T_n}
= \frac{e^{\gamma (\alpha+ w)} }{\sigma (\alpha+ w)}
\frac{\sigma((\alpha+ w)(n+1) + \beta)}{\sigma((\alpha+ w)   n +
\beta)}   \frac{[n]!^2}{[-\hat \alpha]_n [\hat \alpha+2]_n}.
\label{h1_expl} \end{gather}

Formulas \eqref{h_expl} and \eqref{h1_expl} allow us
to f\/ind explicit expressions for the recurrence coef\/f\/i\-cients~$b_n$,~$d_n$.

Indeed, from \eqref{d_n} and \eqref{b_n} we have \begin{gather} d_n
=\frac{h_n^{(1)}}{h_n} = e^{\gamma w}   \frac{[\hat
\alpha-n][\hat \beta+(\hat \alpha+1)(n+1)][\hat \beta+ \hat \alpha
n]}{[\hat \alpha+(n+1)][\hat \beta+(\hat \alpha+1) n][\hat \beta+
\hat \alpha (n+1)]} \label{d_hh} \end{gather} and \begin{gather} b_n = -
\frac{h_n^{(1)}}{h_{n-1}} = -e^{\gamma w}  \frac{[n]^2[\hat
\beta+(\hat \alpha+1) (n+1)][\hat \beta+ \hat \alpha(n-1)]}{[\hat
\beta+ (\hat \alpha+1) n][\hat \beta+ \hat \alpha n ][\hat \alpha
+ n][\hat \alpha + n+1]}. \label{b_hh} \end{gather} We thus obtained a new
explicit example of the Laurent biorthogonal polynomials which
have both explicit expression in terms of the elliptic
hypergeometric function~\eqref{hyp_LP} and explicit recurrence
coef\/f\/icients~\eqref{d_hh},~\eqref{b_hh}.

As a by-product, we have also obtained a new explicit solution of
the discrete-time relativistic Toda chain or, equivalently, a new
explicit solution of the two-point $QD$-algorithm. Indeed, the
recurrence coef\/f\/icients $b_n$, $d_n$ given by~\eqref{d_hh} and~\eqref{b_hh} provide an explicit elliptic solution of the two-point
$QD$-algorithm~\eqref{rel_QD} with $t=\alpha$, $h=w$. In turn, using
correspondence \eqref{W_P_cor} we can obtain elliptic solution of the
ordinary $QD$-algorithm~\eqref{QDR}, or equivalently, the
discrete-time Toda chain solutions.  As far as we know these
solutions are new.

In order to f\/ind explicit (bi)orthogonality relation for these
polynomials we need f\/irst the explicit Fourier expansion of the elliptic
functions of the second kind. We will do this in the next section.

\section{Fourier series of the elliptic functions of the second kind}

Assume that $f(z)$ is the simplest
elliptic function of the second kind \cite{Akhiezer}  \begin{gather} f(z) =
\kappa  \frac{\sigma(z+ \alpha +
\beta)}{\sigma(z+\alpha)}e^{\gamma z} \label{f_pseudo} \end{gather} with some
complex parameters $\kappa$, $\beta$, $\alpha$, $\gamma$. The function
$f(z)$ is quasi-periodic with respect to periods $2 \omega_1$, $2
\omega_3$: \begin{gather} f(z+2 \omega_1) = \mu_1 f(z), \qquad f(z+2 \omega_3)
= \mu_3 f(z), \label{quasi_f} \end{gather} where $\mu_1 = e^{2\eta_1 \beta +
2\omega_1 \gamma}$, $\mu_3 = e^{2\eta_3 \beta + 2\omega_3
\gamma}$. We demand that function $f(z)$ be {\it purely periodic}
with respect to the (real) period $2\omega_1 j$: \begin{gather*} f(z +
2\omega_1 j) = f(z), 
\end{gather*} where $j=1,2, \dots$ is an
arbitrary positive integer. This leads to the condition $\mu_1^j
=1$ or \begin{gather} j(\omega_1 \gamma + \eta_1 \beta) = i \pi m,
\label{per_j_cond} \end{gather} where $m =0, \pm 1, \pm 2, \dots$. Note that
for $j>1$ we should avoid the values $m=0, \pm j, \pm 2j, \dots$
because they correspond to pure $2 \omega_1$-periodicity. Of
course, it is assumed that $m$ and $j$ are coprime, i.e.\ $\mu_1$
is a primitive root of the unity of order $j$: \begin{gather*} \mu_1 =
e^{\frac{2\pi m i}{j}}. 
\end{gather*}

Moreover, we assume that $\alpha= -\alpha_0 - i \alpha_1$, where
both parameters $\alpha_{0,1}$ are real and are restricted by
conditions \begin{gather} 0 \le  \alpha_0 < 2 \omega_1, \qquad 0 < \alpha_1 <
2 |\omega_3|. \label{al_bt_cond} \end{gather} Conditions \eqref{al_bt_cond} mean
that the parameter $-\alpha$ lies within the fundamental
parallelogram (i.e.\ rectangle in our case). If $\alpha$ takes
values beyond this parallelogram, it is possible to reduce it to
canonical choice \eqref{al_bt_cond} using shifts by periods $2
\omega_1$, $2 \omega_3$. Due to quasiperiodicity property of the
function $f(z)$ this will lead only to redef\/ining of the parameter
$\gamma$. Moreover we assume that the imaginary part $-\alpha_1$
of $\alpha$ is nonzero. This assumption is very natural if we
would like to avoid singularities of the function $f(z)$ on whole
real axis. Equivalently, one can present $\alpha$ in the form \begin{gather}
\alpha= -\alpha_0 - 2 \nu \omega_3, \label{alpha_om}
\end{gather} where
$0<\nu<1$ is a f\/ixed parameter which describes the relative value
of the imaginary part $\alpha_1 = - 2i \nu \omega_3$ with respect
to the imaginary period $2 \omega_3$.

Thus we have the function $f(z)$ which is periodic and bounded on
the whole real axis. It is possible therefore to present $f(z)$ in
terms of the Fourier series \begin{gather} f(z) = \sum_{n=-\infty}^{\infty}
A_n \exp\left(\frac{\pi i n z}{j\omega_1}\right).
\label{phi_Fourier} \end{gather} Our problem now is to calculate the Fourier
coef\/f\/icients $A_n$.

By def\/inition, \begin{gather} A_n = \frac{1}{T} \int_{0}^{T} f(z)
\exp\left(\frac{-2\pi i n z}{T}\right) dz, \qquad T= 2 j\omega_1
\label{A_n_def} \end{gather} (the integral is well def\/ined because by our
assumptions the function $f(z)$ has no singularities on the real
axis).

In order to calculate the integral in \eqref{A_n_def} we exploit
standard method of contour integration (see, e.g.,~\cite{Akhiezer}
for calculation of the Fourier expansion for Jacobi elliptic
functions). Choose the contour $\Gamma$ as the rectangle with
vertices $(0$, $2 j\omega_1$, $2j\omega_1+2 \omega_3$, $2\omega_3)$
(i.e.\ the horizontal length is~$2j \omega_1$ and vertical length
$2 |\omega_3|$).

We have (the contour is traversed counterclockwise)
\[
\int_{\Gamma}  f(z)/T \: \exp\left(\frac{-2\pi i n z}{T}\right) dz
= \int_{1} + \int_2 + \int_3 + \int_4,
\]
where $\int_1$, $\int_3$ correspond to horizontal sides of the
rectangle, and integrals $\int_2$, $\int_4$ correspond to vertical
sides.

Due to periodicity property $f(z+ 2 j \omega_1) = f(z)$ we have
$\int_2 + \int_4=0$. For the two remaining horizontal integrals we
have
\[
\int_1 = A_n
\]
and
\[
\int_3= -\int_{2 \omega_3}^{2\omega_3 + T}f(z)/T
\exp\left(\frac{-2\pi i n z}{T}\right) dz .
\]
Making the shift $z \to z + 2\omega_3$ and using quasi-periodic
property \eqref{quasi_f} we have
\[
\int_3=-\mu_3 \exp\left(-\frac{4 \pi i \omega_3 n}{T} \right)
\int_1
\]
and thus \begin{gather} \int_{\Gamma} f(z)/T   \exp\left(\frac{-2\pi i n
z}{T}\right) dz = \left(1- \mu_3 \exp\left(-\frac{4 \pi i \omega_3
n}{T} \right) \right) A_n . \label{int_G} \end{gather} Hence, in order to
calculate the Fourier coef\/f\/icient $A_n$ we need to calculate the
contour integral in l.h.s.\ of~\eqref{int_G}. This can be done by
standard methods of residue theory.

Indeed, inside the contour $\Gamma$ the function $f(z)
\exp\left(\frac{-2\pi i n z}{T}\right)$ has only $j$ simple poles
located at points
\[
z_s= \alpha_0 +i \alpha_1 + 2 s \omega_1, \qquad s=0,1,\dots, j-1 .
\]
At $z_0=-\alpha= \alpha_0 + i \alpha_1$ the function $f(z)$ has
the residue
\[
r=\kappa   e^{-\gamma \alpha} \sigma(\beta).
\] At $z_s$ the
function $f(z)$ has the residue
\[
r_s=\mu_1^s r .
\]
Hence we have that the residue $R_n$ of the function $f(z)/T
\exp\left(\frac{a z}{T}\right)$ inside the rectangle $\Gamma$ will
be
\[
R_n = \frac{r e^{-\chi \alpha}}{T} \sum_{s=0}^{j-1} \mu_1^s
e^{\chi Ts} = \frac{r e^{-\chi \alpha}}{T} \big(1+ q+ q^2 + \cdots +
q^{j-1}\big),
\]
where
\[
\chi= -\frac{\pi i n}{j\omega_1}, \qquad q= \mu_1 e^{\chi T} =
\exp\left(\frac{2 \pi i (m-n)}{j}\right).
\]
If $n \ne m \, \mbox{mod}\, j$ then $R_n=0$. Nonzero value of the
residue will be only for $n =m+ jt$, $t=0, \pm 1, \pm 2,\dots$.
In this case \begin{gather*} R_n = \frac{jr e^{-\chi \alpha}}{T} = \frac{
\kappa \sigma(\beta)}{2 \omega_1} \exp\left(\frac{\alpha \beta
\eta_1}{\omega_1}\right)
  \exp\left(\frac{i \pi \alpha t }{\omega_1}\right). 
  \end{gather*}
Comparing with \eqref{int_G} we get \begin{gather*} A_n = \frac{2 \pi i R_n}{1 -
\mu_3 \exp\left(-\frac{2\pi i \omega_3 n}{j\omega_1}\right)},
\qquad n= m , m \pm j, m \pm 2j ,\dots 
\end{gather*} and
\[
A_n = 0, \qquad \mbox{if} \quad n \ne m \quad \mbox{mod} \; j .
\]
We can simplify this expression using the Legendre identity
\cite{Akhiezer}
\[
\eta_1\omega_3 - \eta_3 \omega_1= \frac{i \pi}{2}
\]
which is valid if ${\rm Im}(\omega_3/\omega_1)>0$. Also we use the
notation \cite{Akhiezer}
\[
h = \exp\left(\frac{i \pi \omega_3}{\omega_1}  \right) .
\]
In our case when $\omega_1 >0$, $i\omega_3<0$ we have that $0<h<1$
(this is so-called normal case for the elliptic function
\cite{Akhiezer}).

We then have
\[
\mu_3 = h ^{2m/j} e^{-\frac{i \pi \beta}{\omega_1}}
\]
and
\[
R_n = R_0  \exp\left(-\frac{i \pi \alpha_0 (n-m)}{j
\omega_1}\right)   h^{\frac{2\nu(m-n)}{j}},
\]
where
\[
R_0 = \frac{ \kappa \sigma(\beta)}{2 \omega_1}
\exp\left(\frac{\alpha \beta \eta_1}{\omega_1}\right)
\]
and we took into account relation \eqref{alpha_om}.

Thus for $n=m+jk$, $k=0, \pm 1, \pm 2, \dots$ we have \begin{gather} A_n=
\frac{2 \pi i R_0  \exp\left(-\frac{i \pi \alpha_0 k}{
\omega_1}\right)   h^{-2\nu k}}{1 - e^{-\frac{i \pi
\beta}{\omega_1}} h^{-2k}} \label{A_n_simpl} \end{gather} and $A_n=0$ if $n
\ne m$ mod\,$j$.

Recall that $j$ is a f\/ixed positive integer -- the order of the
root of unity $\mu_1$, while $m$ is a f\/ixed nonnegative integer
(lesser than $j$) coprime with~$j$. Thus for large $j$ the nonzero
coef\/f\/icients $A_n$ are more rare then for small $j$.

There are two important simplest cases:

(i) if $j=1$ and $m=0$. This case corresponds to the period $2
\omega_1$. Then the Fourier coef\/f\/icients $A_n$ are nonzero for all
$n=0, \pm 1, \pm 2, \dots$ and we have \begin{gather*} A_n =\frac{2 \pi i R_0
  \exp\left(-\frac{i \pi \alpha_0 n}{ \omega_1}\right)   h^{-2n
\nu}}{1 - e^{-\frac{i \pi \beta}{\omega_1}} h^{-2 n}}. 
\end{gather*}

(ii) if $j=2$ and $m=1$. This case corresponds to the period $4
\omega_1$. In this case all even Fourier coef\/f\/icients are zero
$A_{2n}=0$ and for the odd Fourier coef\/f\/icients we have \begin{gather*}
A_{2n+1} =\frac{2 \pi i R_0   \exp\left(-\frac{i \pi \alpha_0 n}{
\omega_1}\right)   h^{-2 n \nu}}{1 - e^{-\frac{i \pi
\beta}{\omega_1}} h^{-2 n}}. 
\end{gather*}

Note that in all cases the Fourier series \eqref{phi_Fourier}
converges inside the strip $-v_1<{\rm Im}(z)<v_2$, where
\[
v_1 = 2|\omega_3| (1-\nu), \qquad v_2 = 2|\omega_3| \nu .
\]
This results follows from standard theorems concerning asymptotic
behavior of the Fourier coef\/f\/icients $A_n$ and $A_{-n}$ for $ n
\to \infty$ \cite{Akhiezer}. The parameters $v_1$, $v_2$ are
positive as follows from the inequality $0<\nu<1$. These
conditions are very natural because the boundary lines ${\rm Im}(z) = 2
\nu |\omega_3|$ and ${\rm Im}(z) = 2 (\nu-1) |\omega_3|$ of the strip
pass through the poles of the function $\phi(z)$. Note that for
$\nu =1/2$ (i.e.\ when the pole of the function $\phi(z)$ lies on
the horizontal line ${\rm Im}(z) = |\omega_3|$) we have the strip
symmetric with respect to the real line: $|{\rm Im}(z)|< |\omega_3|$.
The latter case correspond, e.g., to the Jacobi elliptic functions
$\sn(z;k)$, $\cn(z;k)$, $\dn(z;k)$~\cite{Akhiezer}.

\section{Explicit biorthogonality relation}

In this section we obtain explicit
biorthogonality property of the obtained Laurent biorthogonal
polynomials.

To do this we need to f\/ind explicit realization of the moments
$c_n$ given by formula~\eqref{mom_T}. We note that
\[
c_n = f(wn),
\]
where $f(z)$ is the elliptic function of the second kind
\eqref{f_pseudo} (in our case $\kappa=1/\sigma(\beta)$ but the
constant $\kappa$ does not play any role in formulas for the
polynomials $P_n(z)$ and their recurrence coef\/f\/icients $b_n$,
$d_n$).

Assume f\/irst that the parameter $\gamma$ is chosen to provide the
periodicity of the function $f(z)$ with period $2 \omega_1 j$, $
j=1,2,\dots$. Then we have the Fourier expansion \eqref{phi_Fourier}
from which one obtains \begin{gather} c_n = \sum_{s=-\infty}^{\infty} A_s
\exp\left(\frac{i \pi s w n}{j \omega_1}\right) =
\sum_{s=-\infty}^{\infty} A_s z_s^n,  \label{c_n_A_s} \end{gather} where \begin{gather}
z_s = \exp\left(\frac{i \pi s w}{j \omega_1} \right), \qquad s=0,
\pm 1, \pm 2, \dots \label{z_s} \end{gather} is an inf\/inite set of points
belonging to the unit circle $|z_s|=1$. Due to condition
\eqref{incomp_w} we have that all these points are distinct $z_s \ne
z_t$ if $t \ne s$ and hence they are dense on the unit circle.

From \eqref{c_n_A_s} it follows that the moments $c_n$ are
expressible in terms of the Lebesgue integral \begin{gather*} c_n
=\frac{1}{2\pi} \int_{0}^{2\pi} e^{i \theta n} d \mu(\theta)
\end{gather*} over the unit circle $|z|=1$, where
$\mu(\theta)$ is a (complex) function of bounded variation on the
interval $[0, 2 \pi]$ consisting only from discrete jumps $A_s$
localized in the points $\theta_s$ given by \eqref{z_s}.

Thus we found explicit realization of the moments $c_n$ and hence
we immediately obtain biorthogonality relation for our Laurent
biorthogonal polynomials \begin{gather} \sum_{s=-\infty}^{\infty} A_s
P_n(z_s) Q_m(1/z_s)  = h_n \delta_{nm}, \label{biort_expl} \end{gather} where
$Q_n(z)$ are biorthogonal partners \eqref{deterQ} with respect to
polynomials $P_n(z)$. The Fourier coef\/f\/icients $A_s$ play the role
of discrete weights in this biorthogonality relation. Hence we
have obtained

\begin{proposition}
In the periodic case $f(z+ 2 \omega_1 j) = f(z)$ the elliptic
polynomials \eqref{hyp_LP} $P_n(z)$ are biorthogonal \eqref{biort_expl}
on the unit circle $|z|=1$ with respect to a dense point measure
with weights~$A_s$ given by expression \eqref{A_n_simpl}.
\end{proposition}

Note that the biorthogonal partners $Q_n(z)$ in our case can be
found explicitly in terms of the elliptic hypergeometric function.
Indeed, from \eqref{deterQ} we see that the polynomials $Q_n(z)$ are
Laurent biorthogonal polynomials corresponding to the ``ref\/lected''
moments $\t c_n = c_{-n}$. From explicit expression \eqref{mom_T} it
follows that the moments $c_{-n}$ are obtained from the moments
$c_n$ by ref\/lection of the parameters $\alpha \to - \alpha$, $\beta
\to -\beta$, $\gamma \to - \gamma$, whereas the parameter $w$
remains unchanged (under such procedure we obtain the moments
$-c_{-n}$ but any constant common factor in front of moments leads
to the same polynomials $Q_n(z)$). Hence we can obtain expression
for the polynomials $Q_n(z)$ from the expression~\eqref{hyp_LP} for
polynomials $P_n(z)$ by ref\/lection of parameters $\alpha$, $\beta$,
$\gamma$: \begin{gather*} Q_n(z) = \t B_n \,  {_{3}}E_2 \left( {-n , 1- \hat
\alpha, (\hat \alpha-1) n + \hat \beta +1 \atop 1-n - \hat \alpha,
(\hat \alpha-1)n + \hat \beta}; ze^{\gamma   w} \right),
\end{gather*} where the coef\/f\/icient $\t B_n$ is obtained from
corresponding coef\/f\/icient $B_n$ \eqref{B-cf} by the same ref\/lection
of the parameters $\alpha$, $\beta$, $\gamma$.

Thus both polynomials $P_n(z)$ and their biorthogonal partners
$Q_n(z)$ have similar expressions in terms of elliptic
hypergeometric function.

So far, we assumed that the function $f(z)$ is periodic with the
period $2 \omega_1 j$. This assumption means that the parameter
$\gamma$ should satisfy condition \eqref{per_j_cond}. Parameters
$\alpha$ and $\beta$ are assumed to be arbitrary (with the only
condition \eqref{al_bt_cond}). What happens if the function $f(z)$ is
not periodic, i.e.\ if the parameter $\gamma$ is arbitrary? It
appears that this general case can be easily reduced to the
already considered. Indeed, assume that we change the parameter
$\gamma$, i.e.\ assume that the parameters $\alpha$ and  $\beta$
remain the same but $\t \gamma = \gamma + \chi$, where $\chi$ is
an arbitrary complex parameter. Then it is easily seen from
explicit expression \eqref{hyp_LP} that the new Laurent biorthogonal
polynomials $\t P_n(z)$ are obtained by simple rescaling of the
argument: \begin{gather*} \t P_n(z) = q^n   P_n(z/q), \label{resc_P}
\end{gather*} where
$q = e^{w \chi}$. This corresponds to transformation of the
moments $\t c_n = \epsilon q^n c_n$ as seen directly from
\eqref{mom_T} (the common constant $\epsilon = e^{\alpha w}$ is
inessential and can be put equal to~1).

Assume that we choose the parameter $\chi$ such that the new
function
\[
\t f(z) = \kappa\: \frac{\sigma(z+ \alpha +
\beta)}{\sigma(z+\alpha)}e^{\t \gamma z}
\]
will be periodic with the period $2 \omega_1 j$. This means that
the parameter $\chi$ should be chosen from condition (see
\eqref{per_j_cond}) \begin{gather} \omega_1 (\gamma + \chi) + \eta_1 \beta = i
\pi m/j, \label{cond_chi} \end{gather} where $m$ is co-prime with $j$.

Then the new polynomials $\t P_n(z)$ will be biorthogonal on the
unit circle according to above obtained proposition: \begin{gather}
\sum_{s=-\infty}^{\infty} A_s \t P_n(z_s) \t Q_m(1/z_s) =  h_n
\delta_{nm}, \label{t_P_biort} \end{gather} where the spectral points $z_s$
on the unit circle are given by \eqref{z_s} and the weights $A_s$ by
\eqref{A_n_simpl}. Note that the normalization coef\/f\/icients $h_n$
remain unchanged under the rescaling transform as seen from
\eqref{hbd}, i.e.\ $\t h_n = h_n$.

Taking into account that $\t Q_n(z) = q^{-n} Q_n(z)$ (see
\eqref{res_Q}) we obtain from \eqref{t_P_biort} the biorthogonal
relation \begin{gather} \sum_{s=-\infty}^{\infty} A_s  P_n(z_s/q)  Q_m(q/z_s)
= h_n \delta_{nm}. \label{biort_nunit} \end{gather} Relation \eqref{biort_nunit}
means that for generic values of $\gamma$ polynomials $P_n(z)$ and
$Q_n(z)$ are biorthogonal on the non-unit circle $|z|=1/|q|$ with
respect to the same dense point measure.

It is interesting to note that for every integer $j=1,2,\dots$
(i.e.\ for every period $T=2 \omega_1 j$) we can construct
corresponding circle providing biorthogonality relation
\eqref{biort_nunit}. Thus there exist inf\/initely many orthogonality
circles for dif\/ferent values of the integer parameter~$j$.

For the radius $r$ of the circle of biorthogonality we have from
\eqref{cond_chi} (recall that we assume parameter $w$ to be real) \begin{gather*}
r = 1/|q| = \big| e^{\frac{\eta_1\beta w}{\omega_1}}\big|
\left| e^{w \gamma}\right|. 
\end{gather*}

\section{Positivity of the measure and polynomials orthogonal\\ on the unit circle}

Return to the case when the function
$f(z)$ is periodic with the period $2 \omega_1 j$ and consider an
important special case when all the Fourier coef\/f\/icients of the
function $f(z)$ are nonnegative $A_n \ge 0 $. In this case all
spectral points $z_s$ belong to the unit circle $|z_s|=1$ and the
measure on the unit circle is a positive nondecreasing function.

We have $0<h<1$. Thus for $n \to - \infty$ we have
\[
A_n = 2 \pi i R_0 e^{\frac{-i\pi \alpha_0 n }{\omega_1}}
h^{-2n\nu} .
\]
It is seen that for positivity of $A_n$ one should have $2 \pi i
R_0=\kappa_0$, where $\kappa_0$ is a positive parameter, and for
the real part of $\alpha$ we have the condition  \begin{gather} \alpha_0 =
2J_0\omega_1, \qquad J_0 =0, \pm 1, \pm 2, \dots. \label{cond_alpha} \end{gather}
Now for for $n \to \infty$ we have
\[
A_n= -\kappa_0  e^{\frac{i \pi \beta}{\omega_1}}   h^{2(1-\nu)k}
.
\]
In this case we should have necessarily \begin{gather} {\rm Re}(\beta) = (2J_1+1)
\omega_1, \qquad J_1=0, \pm 1, \pm 2, \dots . \label{cond_beta} \end{gather}
It is easily seen that conditions \eqref{cond_alpha} and
\eqref{cond_beta} are also suf\/f\/icient and so we have the

\begin{proposition}
The Fourier coefficients are positive $($up to inessential common
factor$)$ if and only if the real parts of parameters $\alpha$,
$\beta$ satisfy conditions \eqref{cond_alpha} and \eqref{cond_beta}. In
this case the expression for the Fourier coefficients can be
presented in the form \begin{gather} A_n = \kappa_0  \frac{h^{-2\nu k}}{1+
\kappa_1 h^{-2k}},\qquad n= m + jk, \quad k=0, \pm 1, \pm 2, \dots,
\quad 0<\nu<1, \label{pos_A} \end{gather} and $A_n=0$ if $n \ne m$ ${\rm mod}\,(j)$, where $\kappa_1 =e^{\frac{\pi\, {\rm Im}(\beta)}{\omega_1}}$ is
a positive parameter $($as usual by ${\rm Im}(\beta)$ we denote the
imaginary part of $\beta)$.
\end{proposition}

In this case we have positive dense point measure on the unit
circle. It is well known that when the measure $d \sigma$ is
positive on the unit circle then biorthogonal polynomials become
the orthogonal polynomials on the unit circle \cite{Sz,Ger, Simon}. In this case the moments $c_n$ satisfy
the restriction \begin{gather*} c_{-n} = \bar c_n 
\end{gather*} and
moreover all the Toeplitz determinants are positive \begin{gather*} \Delta_n
>0, \qquad n=1,2,\dots . 
\end{gather*} The property $c_{-n} = \bar c_n$
can be verif\/ied directly from the def\/inition \eqref{mom_T} if the
parameters $\alpha$, $\beta$ satisfy conditions: \begin{gather} \alpha= 2J_1
\omega_1 -  2\nu \omega_3, \qquad \beta= (2J_1+1)\omega_1 + i
\beta_1 \label{al_beta_con} \end{gather} (here $\beta_1$ is an arbitrary real
parameter).

In this special case the obtained polynomials satisfy the
Szeg\H{o} recurrence relation~\eqref{Sz_rec}. The ref\/lection
parameters $a_n$ are calculated as $a_n = -P_{n+1}(0)$ and using
already found explicit formula \eqref{hyp_LP} for polynomials
$P_n(z)$ we have $a_n = -B_{n+1}$, where $B_n$ is given by
\eqref{B-cf} (with~$\alpha$,~$\beta$ satisfying restrictions
\eqref{al_beta_con}). From general theory it follows that in this
case the ref\/lection parameters should satisfy the restriction
$|a_n|<1$. This property is not obvious from explicit expression
for $a_n$ in terms of elliptic Pochhammer symbols.

If, in addition to positivity of $A_n$,  we demand that the
discrete measure should be {\it symmetric} with respect to the
real axis we then obtain the condition $A_{-n}=A_n$ for all $n=0,
1,2, \dots$. It is easily verif\/ied from explicit expression
\eqref{pos_A} that this is possible only for $j=1$ and $j=2$. In the
f\/irst case, when $j=1$ the period $T=2 \omega_1$ and necessarily
$\nu = 1/2$ and $\kappa_1=1$, so that \begin{gather} A_n =
\frac{\kappa_0}{h^n + h^{-n}}. \label{A_n_dn} \end{gather} But the Fourier
coef\/f\/icients with expression \eqref{A_n_dn} correspond to the Jacobi
elliptic function $\dn(z;k)$ \cite{WW}. In this case the moments
are $c_n = \dn(wn;k)$ and indeed satisfy the property
$c_{-n}=c_n$; the ref\/lection parameters are very simple: $a_n=
\dn(w(n+1);k)$ for the even $n$ and $a_n=-\cn(w(n+1);k)$ for the
odd $n$.

In the second case, i.e.\ when $j=2$ we have the period $T=4
\omega_1$ and necessarily $\nu=1/2$ and $\kappa_1=h^{-1}$, so that
\begin{gather*} A_{2n+1} = \frac{\kappa_0}{h^{-n-1/2} + h^{n+1/2}}.
\end{gather*} These Fourier coef\/f\/icients correspond to the
Jacobi elliptic function $\cn(z;k)$ \cite{WW}. Again the moments
$c_n=\cn(wn;k)$ satisfy the desired property $c_{-n}=c_n$ and we
have the polynomials orthogonal on the  unit circle with simple
ref\/lection parameters: $a_n=\cn(w(n+1);k)$ for the even $n$ and
$a_n = -\dn(w(n+1);k)$ for the odd $n$.

These two explicit cases of OPUC with dense point measure were
f\/irst considered in \cite{Zhe_cndn}. Now we see that there exists
much wider class of explicit elliptic OPUC with positive dense
measure on the unit circle. This class of OPUC contains
essentially 3 arbitrary continuous parameters: $w$, ${\rm Im}(\alpha) =
-\alpha_1$, ${\rm Im}(\beta)$. We thus have two additional parameters
with respect to the only parameter $w$ in \cite{Zhe_cndn}. Note
however, that if one demands that OPUC were {\it real} (i.e.\ they
have real ref\/lection parameters $a_n$ and moments $c_n$) then
nothing more general than ``\cn-'' and ``\dn-''polynomials considered
in \cite{Zhe_cndn} appear.

\section[''Classical'' property of LBP]{``Classical'' property of LBP}

\setcounter{equation}{0} Assume that $P_n(z)$ are arbitrary
Laurent polynomials satisfying 3-term RR \begin{gather*} P_{n+1}(z) + d_n
P_n(z) = z(P_n(z) + b_n P_{n-1}(z)) 
\end{gather*} with some
recurrence coef\/f\/icients $b_n$, $d_n$.

For any sequence $\mu_n$, $n=0,1,2,\dots$ of complex numbers such
that $\mu_0=0$ we def\/ine the linear operator $\cal D$ which acts
on the space of polynomials in the argument $z$ by the rule \begin{gather*}
{\cal D} z^n = \mu_n z^{n-1}.  
\end{gather*} Then it is clear
that the operator $\cal D$ sends any polynomial of degree $n$ to a
polynomial of degree $n-1$ and moreover ${\cal D} \{1\} =0$. In
this sense the operator $\cal D$ can be called as a generalized
derivative operator. If $\mu_n=n$ then ${\cal D} = \partial_z $
coincides with the ordinary derivative operator with respect to
the variable $z$.

We say that LBP $P_n(z)$ are $\cal D$-classical if \begin{gather} {\cal D}
P_n(z) = \mu_n \t P_{n-1}(z), \label{D_clas} \end{gather} where $\t P_n(z)$
is another set of LBP satisfying the recurrence relation \begin{gather*} \t
P_{n+1}(z) + \t d_n \t P_n(z) = z(\t P_n(z) + \t b_n \t
P_{n-1}(z)) 
\end{gather*} with some coef\/f\/icients $\t d_n$, $\t b_n$
and initial conditions $\t P_0=1$, $\t P_1 = z- \t d_0$.

In \cite{Zh_LBP} we considered the case of the ordinary classical
LBP (i.e.\ with respect to the opera\-tor~$\partial_z$) and derived
necessary and suf\/f\/icient conditions for existence of such
polynomials. It appears that there exists many dif\/ferent types of
such classical LBP. The simplest ones are the LBP constructed by
Hendriksen and van~Rossum \cite{HR}. The latter have explicit
expression in terms of the Gauss hypergeometric function.

Now return to our elliptic LBP $P_n(z)=P_n(z;\alpha, \beta,
\gamma,w)$ (we indicate dependence on parameters $\alpha$, $\beta$,
$\gamma$, $w$ for convenience) and consider the operator $\cal D$
with $\mu_n$ def\/ined as \begin{gather} \mu_n = \frac{\sigma(wn)}{\sigma(wn +
\alpha)}. \label{ell_mu} \end{gather} Then from explicit representation
\eqref{hyp_LP} it is elementary to verify that the operator $\cal D$
transforms these Laurent biorthogonal polynomials to the same
family but with the sole parameter $\beta$ changed: \begin{gather*} {\cal D}
P_n(z;\alpha, \beta, \gamma,w) = \mu_n   P_{n-1}(z; \alpha,
\beta- \alpha, \gamma,w) 
\end{gather*} which means that our
elliptic polynomials $P_n(z;\alpha, \beta, \gamma,w)$ are indeed
``classical'' polynomials with respect to the operator $\cal D$.

In particular, the choice $\alpha=-\omega_3$, $\beta=\omega_1$
corresponds to the OPUC $\dn$-polynomials \cite{Zhe_cndn}.  Under
action of the operator~$\cal D$ we obtain polynomials with $
\alpha=-\omega_3$, $\beta=\omega_1+\omega_3=-\omega_2$. These
polynomials correspond to the OPUC $\cn$-polynomials
\cite{Zhe_cndn}. Vice versa, action of the opera\-tor~$\cal D$ on
the $\cn$-polynomials return them to $\dn$-polynomials. In this
case the coef\/f\/icient~$\mu_n$ is proportional to the Jacobi
$\sn$-function: $\mu_n = {\rm const}\cdot \sn(wn;k)$~\cite{Zhe_cndn}.

One can repeat action of the operator $\cal D$. This leads to a
chain of corresponding transformations of polynomials $P_n(z)$:
\begin{gather*} {\cal D}^m   P_n(z;\alpha, \beta, \gamma,w) = \mu_n \mu_{n-1}
\cdots \mu_{n-m+1} P_{n-m}(z; \alpha, \beta-m \alpha, \gamma, w).
\end{gather*}

One can consider the ``$\mu$-exponential'' function $E_{\mu}(x)$
which is a formal solution of the opera\-tor equation \begin{gather} {\cal D}
E_{\mu}(x) = E_{\mu}(x). \label{E_mu_eq} \end{gather}  Clearly we have a
solution of the opera\-tor equation \eqref{E_mu_eq} in terms of the
formal series \begin{gather*} E_{\mu}(x) = \sum_{s=0}^{\infty}
\frac{x^s}{\mu_1 \mu_2 \cdots \mu_s}. 
\end{gather*} In case of
the elliptic $\mu_n$ \eqref{ell_mu} we have \begin{gather} E_{\mu}(x) =
\sum_{s=0}^{\infty} \frac{[\alpha+1]_s}{[s]!}    x^s .
\label{E_ell}
\end{gather} Function \eqref{E_ell} is closely related with
so-called ``theta analogue'' of the exponential function proposed by
Spiridonov in \cite{Spi_Bailey}\footnote{The authors are indebted to
V.~Spiridonov for drawing their attention to this result.}.

Obviously for an arbitrary complex parameter $\gamma$ we have
\[
{\cal D} E_{\mu}(\gamma x) = \gamma E_{\mu}(\gamma x) .
\]

We can also introduce ``even'' and ``odd'' $\mu$-exponential functions
which are $\mu$-analogs of the hyperbolic ``cosh'' and ``sinh''
functions
\[
C_{\mu}(x) = \sum_{s=0}^{\infty} \frac{
x^{2s}}{\mu_1 \mu_2 \cdots \mu_{2s}}, \qquad S_{\mu}(x) =
\sum_{s=0}^{\infty} \frac{ x^{2s+1}}{\mu_1 \mu_2 \cdots
\mu_{2s+1}}.
\] We have obvious relations
\[
C_{\mu}(x) = (E_{\mu}(x) + E_{\mu}(-x))/2, \qquad S_{\mu}(x) =
(E_{\mu}(x) - E_{\mu}(-x))/2 .
\]

These functions both have the same property
\[
{\cal D}^2 C_{\mu}(\gamma x) = \gamma^2 C_{\mu}(x), \qquad {\cal
D}^2 S_{\mu}(\gamma x) = \gamma^2 S_{\mu}(x)
\]
with an arbitrary parameter $\gamma$, and hence for arbitrary
parameters $\beta_0$, $\beta_1$ the function $f(x)=\beta_0 C_{\mu}(\gamma
x) + \beta_1 S_{\mu}(\gamma x)$ is a formal solution of the
operator equation
\[
{\cal D}^2 f(x) = \gamma^2 f(x) .
\]
Note the obvious ``intertwining'' property of these functions:
\[
{\cal D} C_{\mu}(\gamma x) = \gamma S_{\mu}(\gamma x), \qquad {\cal
D} S_{\mu}(\gamma x) = \gamma C_{\mu}(\gamma x) .
\]
We can use these properties in order to construct formal
generating functions for the $\cn$ and $\dn$-circle polynomials.

Indeed, let $P_n^{(C)}(z)$ and $P_n^{(D)}$ be $\cn$ and
$\dn$-circle polynomials corresponding to the choices
$\alpha=-\omega_3$, $\beta=-\omega_2$ and $\alpha=-\omega_3$, $
\beta=\omega_1$ respectively. As shown in~\cite{Zhe_cndn} these
polynomials satisfy intertwining properties
\[
{\cal D} P_n^{(C)}(z) = \mu_n   P_{n-1}^{(D)}(z), \qquad {\cal D}
P_n^{(D)}(z) = \mu_n   P_{n-1}^{(C)}(z)
\]
and hence
\[
{\cal D}^2 P_n^{(C)}(z) = \mu_n \mu_{n-1} P_{n-2}^{(C)}(z), \qquad
{\cal D}^2 P_n^{(D)}(z) = \mu_n \mu_{n-1} P_{n-2}^{(D)}(z),
\]
where one can choose $\mu_n = \sn(wn)/\sn(w)$.

Construct the generating functions for the polynomials
$P_n^{(C)}(z)$ and $P_n^{(D)}(z)$ as the formal series \begin{gather}
F^{(C)}(z;t) = \sum_{n=0}^{\infty} \frac{t^n P_n^{(C)}(z)}{\mu_1
\mu_2 \cdots \mu_n}, \qquad F^{(D)}(z;t) = \sum_{n=0}^{\infty}
\frac{t^n P_n^{(D)}(z)}{\mu_1 \mu_2 \cdots \mu_n}. \label{gen_cd} \end{gather}
We have obviously
\[
{\cal D}_z F^{(C)}(z;t) = t F^{(D)}(z;t), \qquad {\cal D}_z
F^{(D)}(z;t) = t F^{(C)}(z;t),
\]
where notation ${\cal D}_z$ means that the operator ${\cal D}$
acts only on the variable $z$. As a consequence
\[
{\cal D}_z^2 F^{(C)}(z;t) = t^2 F^{(C)}(z;t), \qquad {\cal D}_z^2
F^{(D)}(z;t) = t^2 F^{(D)}(z;t) .
\]
This property means that both functions $F^{(C)}(z;t)$ and
$F^{(D)}(z;t)$ can be expressed in terms of ``even'' and ``odd''
$\mu$-exponential functions with respect to the variable $z$:
\begin{gather}
F^{(C)}(z;t) = \xi_0(t) C_{\mu}(zt) + \xi_1(t) S_{\mu}(zt),
\label{gen_cn} \\
 F^{(D)}(z;t) = \eta_0(t) C_{\mu}(zt) +
\eta_1(t) S_{\mu}(zt) \label{gen_dn}
\end{gather} with some functions
$\xi_i(t)$, $\eta_i(t)$, $i=1,2$. For these functions we have
relations $\eta_0(t)=\xi_1(t)$, $\eta_1(t) = \xi_0(t)$ which
follow easily from intertwining relations. Remaining functions
$\xi_0(t)$, $\xi_1(t)$ can be found as follows. Put $z=0$. Then from
def\/inition \eqref{gen_cd} we have \begin{gather*} F^{(C)}(0;t) =
\sum_{n=0}^{\infty} \frac{t^n B_n^{(C)}}{\mu_1 \mu_2 \cdots \mu_n},
\qquad F^{(D)}(0;t) = \sum_{n=0}^{\infty} \frac{t^n
B_n^{(D)}}{\mu_1 \mu_2 \cdots \mu_n},
\end{gather*} where
$B_n^{(C)}$, $B_n^{(D)}$ are corresponding normalization
coef\/f\/icients~\eqref{B-cf}. On the other hand from \eqref{gen_cn} and~\eqref{gen_dn} we have{\samepage
\[
F^{(C)}(0;t) = \xi_0(t), \qquad F^{(D)}(0;t) = \xi_1(t)
\]
because $C_{\mu}(0)=1$, $S_{\mu}(0)=0$.}

We thus have explicit expressions for the functions $\xi_0(t)$,
$\xi_1(t)$ in terms of formal series \begin{gather*} \xi_0(t) =
\sum_{n=0}^{\infty} \frac{t^n B_n^{(C)}}{\mu_1 \mu_2 \cdots \mu_n},
\qquad \xi_1(t) = \sum_{n=0}^{\infty} \frac{t^n B_n^{(D)}}{\mu_1
\mu_2 \cdots \mu_n}.
\end{gather*}

As shown in \cite{Zhe_cndn} the coef\/f\/icients $B_n^{(C)}$ and
$B_n^{(D)}$ are expressible in terms of the Jacobi elliptic
functions as follows: \begin{gather*} B_n^{(C)} = - a_{n-1}^{(C)} = \left\{\begin{array}{ll}
-\cn(wn), & n \  \mbox{even}, \\ \dn(wn), & n \
\mbox{odd} \end{array} \right . 
\end{gather*} and \begin{gather*} B_n^{(D)} = -
a_{n-1}^{(D)} = \left\{ \begin{array}{ll} -\dn(wn), & n \  \mbox{even},    \\
\cn(wn), &  n \  \mbox{odd}, \end{array} \right .   
\end{gather*} where
$a_{n}^{(C)}$ and $a_{n}^{(D)}$ are corresponding ref\/lection
parameters.

We can thus present expressions for $\xi_0(t)$ and  $\xi_1(t)$ in
a more explicit form: \begin{gather*} \xi_0(t) = \sum_{s=0}^{\infty} \frac{
\cn(2ws) t^{2s} \sn^{2s}(w)}{\sn(w) \sn(2w) \cdots \sn(2ws)} -
\sum_{s=0}^{\infty} \frac{ \dn(w(2s+1)) t^{2s+1}
\sn^{2s+1}(w)}{\sn(w) \sn(2w) \cdots \sn(w(2s+1))},\\ 
  \xi_1(t) = \sum_{s=0}^{\infty} \frac{ \dn(2ws) t^{2s}
\sn^{2s}(w)}{\sn(w) \sn(2w) \cdots \sn(2ws)} - \sum_{s=0}^{\infty}
\frac{ \cn(w(2s+1)) t^{2s+1} \sn^{2s+1}(w)}{\sn(w) \sn(2w) \cdots
\sn(w(2s+1))}. 
\end{gather*} Hence we were able to f\/ind
explicitly the generating functions for the OPUC $P_n^{(C)}(z)$
and $P_n^{(D)}(z)$ in terms of $\mu$-exponential functions
\eqref{gen_cn}, \eqref{gen_dn}. The problem of convergence of these
functions needs a separate investigation.

\section{Rational limit of the elliptic functions\\ and corresponding Laurent biorthogonal polynomials}

 In this section we consider the rational
limit of the elliptic functions, when both periods~$2\omega_1$,~$2\omega_2$ tend to inf\/inity. In this case for the Weierstrass
functions we have simple formulas \cite{Akhiezer}
\[
\sigma(z) = z, \qquad \zeta(z) = 1/z, \qquad \wp(z) = 1/z^2 .
\]
Hence for the moments $c_n$ \eqref{mom_T} we have (one can put $w=1$
and $\gamma=0$ without loss of generality) \begin{gather} c_n =
\frac{n+\beta+\alpha}{\beta(n+\alpha)}. \label{rat_mom} \end{gather} In this
case the elliptic Frobenius determinant $\Delta_n$ becomes the
well-known rational Cauchy determinant \cite{Krat} and for
corresponding Laurent biorthogonal polynomials \eqref{PL_det} we have
($\alpha_1 = \alpha+1$) \begin{gather} P_n(z) = B_n  \, {_{3}}F_2 \left( {-n ,
\alpha_1, -\alpha_1 n - \beta +1 \atop \alpha_1 - n, -\alpha_1 n -
\beta}; z \right) \label{rat_LP} \end{gather}  with \begin{gather*} B_n =
\frac{(-\alpha)_n}{(\alpha +1)_n}   \frac{\alpha n + \beta +
n}{\alpha n + \beta},   
\end{gather*} where ${_{3}}F_2(z)$ is
the ordinary hypergeometric function \cite{KS} and
$(a)_n=a(a+1)\cdots (a+n-1)$ is the Pochhammer symbol (shifted
factorial). These polynomials satisfy the 3-term recurrence
relation \eqref{recP} with \begin{gather} d_n =
\frac{(\alpha-n)(\beta+\alpha_1(n+1))(\beta+\alpha
n)}{(\alpha_1+n)(\beta+\alpha_1 n)(\beta + \alpha(n+1))}
\label{d_n_rat} \end{gather} and \begin{gather} b_n= -  \frac{n^2(\beta+ \alpha_1 (n+1))
(\beta+\alpha(n-1))}{(\beta+\alpha_1 n)(\beta+\alpha
n)(\alpha+n+1)(\alpha+n)}. \label{b_n_rat} \end{gather} Note that as in the
elliptic case we should require that \begin{gather*} \alpha \ne 0, \pm 1, \pm
2, \dots 
\end{gather*} because otherwise the normalization
coef\/f\/icients $h_n$ become singular at some $n>0$ which means a
degeneration.

As in the elliptic case the recurrence coef\/f\/icients \eqref{d_n_rat}
and \eqref{b_n_rat} provide an explicit solution of the two-point
$QD$-algorithm \eqref{rel_QD} with $t=\alpha$, $h=1$. Similarly, one can
construct corresponding solutions of the ordinary $QD$-algorithm
\eqref{QDR}.

The biorthogonal partners have the expression \begin{gather*} Q_n(z) = B_n  \,
{_{3}}F_2 \left( {-n , 1-\alpha, (\alpha-1) n + \beta +1 \atop
1-\alpha - n, (\alpha-1) n + \beta}; z \right). 
\end{gather*}

In order to f\/ind the orthogonality measure for these polynomials
we f\/irst note that the moments \eqref{rat_mom} can be rewritten in
the form \begin{gather} c_n = \beta^{-1} + \frac{1}{n+\alpha} = \beta^{-1} +
c_n^{(0)}, \label{mom_HR} \end{gather} where the moments
\[
c_n^{(0)}=1/(n+\alpha) = \lim_{\beta \to \infty} c_n
\]
correspond to a special case of the ``classical'' Laurent
biorthogonal polynomials considered by Hendriksen and Van Rossum
\cite{HR}. The moments $c_n^{(0)}$ correspond to the recurrence
coef\/f\/icients \begin{gather} d_n = \frac{\alpha-n}{\alpha+n+1}, \qquad b_n =
-\frac{n^2}{(\alpha+n)(\alpha+n+1)}. \label{db_HR} \end{gather} The
polynomials $P_n(z)$ are expressed in terms of the Gauss
hypergeometric function \cite{HR} \begin{gather} P_n(z) =
\frac{(-\alpha)_n}{(\alpha+1)_n}  \, {_{2}}F_1 \left( {-n , \alpha
+ 1,  \atop \alpha + 1 - n}; z \right). \label{HR_LP} \end{gather} Formulas
\eqref{HR_LP} and \eqref{db_HR} are easily obtained from \eqref{d_n_rat},
\eqref{b_n_rat} and \eqref{rat_LP} in the limit $\beta \to \infty$. The
biorthogonality property for polynomials \eqref{HR_LP} was found in
\cite{HR} \begin{gather*} \int_{C}P_n(z) Q_m(1/z) z^{\alpha-1}dz = h_{n}
\delta_{nm}, 
\end{gather*} where the contour $C$ is the unit
circle and integrating path starts from $z=1^+$, where
$z^{\alpha-1}=1$ and stops at $z=1^-$, where $z^{\alpha-1} =
e^{2\pi i \alpha}$.

The moments \eqref{mom_HR} dif\/fer from the moments $c_n^{(0)}$ only
by the constant term $\beta^{-1}$. This means that corresponding
orthogonality measure on the unit circle has an additional
concentrated mass at $z=1$. Thus the biorthogonality relation for
polynomials \eqref{rat_LP} looks as \begin{gather*} \int_{C}P_n(z) Q_m(1/z)
w(z)dz = h_{n} \delta_{nm},
\end{gather*} where the weight
function is
\[
w(z) = \frac{z^{\alpha-1}}{e^{2\pi i \alpha}-1} + \beta^{-1}
\delta(z-1).
\] Note that inserting the concentrated mass at point
$z=1$ can be performed by the Geronimus transform of the classical
Hendriksen--van~Rossum polynomials (see \cite{Zh_LBP} for details).

The ``classical'' property \eqref{D_clas} holds for the polynomials
\eqref{rat_LP} with
\[
\mu_n=\frac{n}{n+\alpha} .
\]
In more details
\[
{\cal D} P_n(z; \alpha, \beta) = \mu_n P_{n-1}(z; \alpha,
\beta-\alpha),
\]
where $P_n(z;\alpha,\beta)$ stands for LBP \eqref{rat_LP} with
explicit dependence on the parameters $\alpha$, $\beta$. Note the
operator $\cal D$ in this case does not coincide with the ordinary
derivative operator $\partial_z$. Hence the polynomials
\eqref{rat_LP} provide one of the simplest nontrivial examples of the
``classical'' LBP with respect to a nonclassical ``derivative''
operator $\cal D$.

Consider also a special case when both $\alpha$ and $\beta$ are
purely imaginary parameters: \begin{gather*} \alpha=i s_1, \qquad \beta=i s_2,
\qquad s_1s_2 >0 . 
\end{gather*} Then we can put (obviously
the moments $c_n$ are def\/ined up to a unnecessary common constant
factor) by the expression \begin{gather*} c_n = \frac{n+i(s_1+s_2)}{n+i s_1}
\end{gather*} from which the condition $\bar c_n = c_{-n}$
follows which means that the corresponding polyno\-mials~$P_n(z)\!$
satisfy the Szeg\H{o} recurrence relation \eqref{Sz_rec} and one can
expect that these polynomials are ortho\-go\-nal on the unit circle
with respect to a positive measure. According to the general
theory this occurs if and only if the ref\/lection parameters
$a_n=-P_{n+1}(0)$ satisfy the restriction $|a_n|<1$ for all $n$.

from \eqref{rat_LP} we see that $a_{n-1} =-B_n$. Then
\[
|a_{n-1}|^2= B_n \bar B_n .
\]
After simple calculations we get \begin{gather*} |a_{n-1}|^2 =
\frac{s_1^2}{n^2+s_1^2} \frac{(s_1 n+s_2)^2 + n^2}{(s_1 n+s_2)^2}
= \frac{1+ \xi_n^2}{1+\eta_n^2},  
\end{gather*} where
\[
\eta_n = n/s_1, \qquad \xi_n = n/(n s_1 + s_2) .
\]
It is easily seen that the condition $|a_{n}|<1$ is equivalent to
the condition $s_1s_2>0$. In this case we can present the
orthogonality property in the form \begin{gather*} \int_{\theta=0}^{2\pi}
P_n(e^{i \theta}) \bar P_m(e^{-i \theta}) \rho(\theta) d \theta =
h_n \delta_{nm}, 
\end{gather*} where \begin{gather*} \rho(\theta) =
\frac{s_2}{1-e^{-2\pi s_1}}e^{-s_1 \theta}  + \delta(\theta) .
\end{gather*}
Using correspondence \eqref{W_P_cor} consider the ordinary orthogonal
polynomials $W_n^{(j)}(z)= P_n^{(j+n)}(z)$, where by
$P_n^{(j)}(z)$ we mean polynomials \eqref{rat_LP} obtained under the
substitution $\alpha \to \alpha+j$.

We have \begin{gather} W_n^{(j)}(z) = B_n \, {_{3}}F_2 \left( {-n ,
\alpha+n+j+1, 1- \beta -n(\alpha+j+n+1) \atop \alpha+j+1, - \beta
-n(\alpha+j+n+1)}; z \right), \label{rat_W} \end{gather} with \begin{gather*} B_n =
\frac{(-\alpha-j-n)_n}{(\alpha +j+n+1)_n}   \frac{(\alpha+n+j+1)
n + \beta}{(\alpha+n+j) n + \beta}.  
\end{gather*} These
polynomials are orthogonal with respect to the moment sequence
\[
\tau_n^{(j)} = \beta^{-1} + \frac{1}{n+\alpha+j+1}.
\]
In order to get corresponding orthogonality measure we f\/irst
consider the limit $\beta \to \infty$. In this limit we have the
polynomials \begin{gather*} W_n^{(j)}(z) = B_n \, {_{2}}F_1 \left( {-n ,
\alpha+n+j+1  \atop \alpha+j+1}; z \right) 
\end{gather*}
with \begin{gather*} B_n = \frac{(-\alpha-j-n)_n}{(\alpha +j+n+1)_n}.
\end{gather*} These polynomials coincide with a special
class of Jacobi polynomials which are orthogonal on the interval
$[0,1]$ with the weight function  $w(x) = x^{\alpha+j}$: \begin{gather*}
\int_0^1 x^{\alpha+j} W_n^{(j)}(x) W_m^{(j)}(x) dx = q_n
\delta_{nm} . 
\end{gather*} Indeed, for the moments we have
\[
\int_0^1 x^{\alpha+j} x^n dx = \frac{1}{\alpha+j+n+1} =
\lim_{\beta \to \infty} \tau_n^{(j)}.
\]
Hence, the polynomials $W_n^{(j)}$ with nonzero $\beta$ correspond
to adding a concentrated mass $M=\beta^{-1}$ at the endpoint
$x=1$ of orthogonality interval for the Jacobi polynomials. I.e.\
the weight function for the polynomials $W_n^{(j)}(x)$ with f\/inite
values of $\beta$ has the expression
\[
w(x) = x^{\alpha+j} + \beta^{-1} \delta(x-1).
\]
Such polynomials are called the Krall--Jacobi polynomials (see,
e.g.,~\cite{LL}). These polyno\-mials have a remarkable property:
they are eigenfunctions of a fourth-order dif\/ferential opera\-tor~\mbox{\cite{LL,LL2}}.

Hence we can expect that our Laurent biorthogonal polynomials will
satisfy corresponding fourth-order dif\/ferential equation too. We
f\/irst consider a more elementary case $\beta = \infty$. Then the
Jacobi polynomials $W_n^{(j)}(x)$ satisfy the dif\/ferential
equation \cite{KS} \begin{gather} x(1-x)\frac{d^2 W_n^{(j)}(x)}{dx^2} +
\frac{1}{2} (\alpha+j + (\alpha+j+2)(1-2x) )\frac{d
W_n^{(j)}(x)}{dx}\nonumber\\
\qquad{} + n(n+\alpha+j+1) W_n^{(j)}(x) =0 .
\label{diff_Jac}
\end{gather} In order to return to the Laurent biorthogonal
polynomials we need to shift the parameter $j \to j-n$. Then
equation \eqref{diff_Jac} becomes \begin{gather*} z(1-z)\frac{d^2
P_n^{(j)}(z)}{dz^2} + \frac{1}{2} (\alpha+j-n +
(\alpha+j-n+2)(1-2z) )\frac{d P_n^{(j)}(z)}{dz}\nonumber\\
\qquad{}{} + n(\alpha+j+1)
P_n^{(j)}(z) =0. 
\end{gather*} We can rewrite this equation
in the form of generalized eigenvalue problem \begin{gather} A P_n^{(j)}(z) =
\lambda_n B P_n^{(j)}(z) \label{gevp_rat} \end{gather} for two dif\/ferential
operators
\[
A= z(1-z) \partial_z^2 + (\alpha+j+1 -(\alpha+j+2)z)\partial_z,
\qquad B=(1-z)\partial_z  -\alpha-j-1,
\]
where $\lambda_n=n$ is the generalized eigenvalue.

Consider now general case of f\/inite values of the parameter
$\beta$. The orthogonal polyno\-mials~\eqref{rat_W} satisfy the
dif\/ferential equation of the 4th order \cite{LL, LL2} \begin{gather}
L W_n^{(j)}(x) = \lambda_n W_n^{(j)}(x), \label{LW} \end{gather} where the
dif\/ferential operator $L$ is
\[
L = (x(1-x))^2 \partial_x^4 + x(1-x)(\xi_1 x+ \xi_0) \partial_x^3
+ (1-x)(\eta_1 x + \eta_0) \partial_x^2 + (\zeta_1 x + \zeta_0)
\partial_x
\]
with the coef\/f\/icients
\begin{gather*}
\xi_1= -2(\alpha+j+4), \qquad \xi_0 = 2(\alpha+j+2),\\
\eta_1 =
-(\alpha+j+2)(\alpha+j+7) - 2 \beta, \qquad \eta_0
=(\alpha+j+1)(\alpha+j+2),
\\
\zeta_1=2(\alpha+j+2)(\alpha+\beta+j+1), \qquad \zeta_0 =
-2(\alpha+j+1)(\alpha+\beta+j+2).
\end{gather*}
The eigenvalue is
\begin{gather*}
\lambda_n = n(n-1)(n-2)(n-3) - \xi_1 n(n-1)(n-2) - \eta_1 n(n-1) +
\zeta_1 n
\\
\phantom{\lambda_n}{}=n(n+\alpha+j+1)(n^2 +n(\alpha+j+1) + 2\beta) .
\end{gather*}
If we now substitute $j \to j-n$ into above formulas, we return to
the Laurent biorthogonal polynomials $P_n^{(j)}(z)$ which satisfy
a dif\/ferential equation of the 4th degree which can be presented
in the form \begin{gather} \big(n^2 L_2 + n L_1 + L_0\big) P_n^{(j)}(z) =0,
\label{L_pencil} \end{gather} where $L_0$ is a dif\/ferential operator of the
4th order, $L_1$ of the 3-rd order and  $L_2$ of the second
order. The operators $L_0$, $L_1$, $L_2$ do not depend on the parameter
$n$. We thus see that the Laurent biorthogonal polynomials
$P_n^{(j)}(z)$ given by~\eqref{rat_LP}, satisfy the {\it quadratic}
operator pencil equation~\eqref{L_pencil} (with respect to the
``eigenvalue'' parameter $n$). When $\beta = \infty$ we have
generalized eigenvalue problem \eqref{gevp_rat} which is equivalent
to a {\it linear} operator pencil.

Hence the Laurent biorthogonal polynomials \eqref{L_pencil} can be
considered as biorthogonal ana\-logs of the Krall--Jacobi orthogonal
polynomials. They possess many useful properties of the
Krall--Jacobi polynomials including the 4th order dif\/ferential
equation they satisfy. The main dif\/fe\-ren\-ce, however, is that in
the biorthogonal case we have quadratic operator pencil equation~\eqref{L_pencil} instead of usual eigenvalue problem~\eqref{LW}.

We thus see that many nontrivial properties of the elliptic
Laurent biorthogonal polyno\-mials are manifested already in the
rational limit. Loosely speaking, one can say that elliptic
biorthogonal polynomials \eqref{hyp_LP} are elliptic analogs of the
Krall--Jacobi polynomials.

\subsection*{Acknowledgements}
 The authors thank L.~Golinskii,
A.~Kirillov, C.~Krattenthaler, A.~Magnus, M.~Rahman  and V.~Spiridonov
for discussion. The authors are also indebted to the referees for
careful reading the manuscript and many suggestions leading to
improving of the text.

\pdfbookmark[1]{References}{ref}
\LastPageEnding

\end{document}